\documentclass[12pt,reqno]{amsart}

\usepackage{fullpage,amsfonts,amsmath,amscd,amssymb}

\input{xy}

\xyoption{all}

\theoremstyle{plain}

\newtheorem*{trm}{Theorem}

\newtheorem*{lem}{Lemma}

\newtheorem*{prop}{Proposition}

\newtheorem*{thm}{Theorem}

\newtheorem*{example}{Example}

\newtheorem*{cor}{Corollary}

\theoremstyle{remark}

\newtheorem*{defn}{\bf{Definition}}

\newtheorem*{rem}{\bf{Remark}}

\newcommand{\Z}{\mathbb{Z}}

\newcommand{\R}{\mathbb{R}}

\newcommand{\N}{\mathbb{N}}

\newcommand{\ep}{\epsilon}

\newcommand{\C}{\mathbb{C}}

\newcommand{\s}{\mathcal{S}}

\newcommand{\Co}{\mathcal{O}}

\newcommand{\el}{\mathcal{L}}

\newcommand{\g}{\mathfrak{g}}

\newcommand{\h}{\mathfrak{h}}

\DeclareMathOperator{\gp}{G} \DeclareMathOperator{\ggp}{\hat{G}}

\DeclareMathOperator{\gl}{GL} \DeclareMathOperator{\rep}{Rep}

\DeclareMathOperator{\mat}{Mat} \DeclareMathOperator{\tr}{tr}

\DeclareMathOperator{\edo}{End} \DeclareMathOperator{\lie}{Lie}

\DeclareMathOperator{\sym}{Sym} 

\DeclareMathOperator{\homo}{Hom} 

\DeclareMathOperator{\aut}{Aut} \DeclareMathOperator{\id}{Id}

\DeclareMathOperator{\mx}{Max}

\DeclareMathOperator{\symp}{Symp}

\DeclareMathOperator{\vect}{Vect}

\newcommand{\hmm}[3]{\homo_{#1}(#2,#3)}

\newcommand{\rp}[2]{\rep(#1,#2)}

\newcommand{\ed}[1]{\edo(#1)}

\title{Stratifications of Marsden-Weinstein reductions for representations of quivers and deformations of
symplectic quotient singularities}
\author{Maurizio Martino}

\address{Department of Mathematics, University
of Glasgow, Glasgow, G12 8QW, U.K.}

\email{mma@maths.gla.ac.uk}

\begin{document}
\bibliographystyle{plain}

\begin{abstract}
We investigate the Poisson geometry of the Marsden-Weinstein
reductions of the moment map associated to the cotangent bundle of
the space of representations of a quiver. We show that the
stratification by representation type equals the stratification by
symplectic leaves. The deformed symplectic quotient singularities -
spectra of centres of symplectic reflection algebras associated to a
wreath product - are shown to be isomorphic to reductions of a
certain quiver. This establishes a method to calculate when these
deformations are smooth. Furthermore, the isomorphism identifies
symplectic leaves and so one can give a description of their
symplectic leaves in terms of roots of the quiver.
\end{abstract}

\maketitle

\section{Introduction}

\subsection{} The representation
theory of quivers is a fundamental topic in modern algebra. An
interesting area in this subject has been the study of the geometry
of the moment map, or, equivalently, the representation spaces of
deformed preprojective algebras, with applications to Kleinian
singularities and differential equations. Many interesting results
have been obtained using the combinatorics of root systems for
quivers, see \cite{cb1} for example. We take a different point of
view and examine the symplectic geometry of the moment map. In
particular we describe the symplectic leaves of these
Marsden-Weinstein reductions.

We apply these results to symplectic reflection algebras. These were
introduced by Etingof and Ginzburg in \cite{eg} and are a beautiful
class of algebras which have connections with algebraic geometry,
integrable systems and combinatorics. At parameter $t=0$ the spectra
of the centres of these algebras are Poisson deformations of
symplectic quotient singularities. The symplectic leaves tell us
valuable information about these varieties and also about the
representation theory of the corresponding symplectic reflection
algebras, \cite{brgo}. By using quivers we can calculate, in a large
class of examples, the symplectic leaves in terms of the relatively
well understood theory of roots.

\subsection{} We summarise our results; full details and definitions
are given in later chapters. Let $Q$ be a quiver with vertex set $I$
and set of arrows, $A$, and let $\overline{Q}$ be the double quiver
of $Q$. Let $\lambda\in \C^I$. Then we define the deformed
preprojective algebra
\[ \Pi_{\lambda} = \frac{\C\overline{Q}}{(\sum_{a\in A}[a,a^{\ast}] -
\sum_{i\in I}\lambda_ie_i)}.
\]
The variety
\[ \rep (\Pi_{\lambda}, \alpha) \subset \bigoplus_{a\in \overline{A}} \mat(\alpha_{h(a)} \times
\alpha_{t(a)}, \C) = \rep (\overline{Q}, \alpha)
\]
classifies the representations of $\Pi_{\lambda}$ with dimension
vector $\alpha \in \N^I$ where $\lambda \cdot \alpha = 0$. The group
\[ \gp (\alpha) = (\prod_{i\in I} \mathrm{GL}(\alpha_i,
\C))/\C^{\times}
\]
acts on $\rep (\Pi_{\lambda}, \alpha)$ by basechange. The points of
the algebraic quotient, \begin{equation*} \mathcal{N}(\lambda,
\alpha):= \rep (\Pi_{\lambda}, \alpha)//\gp(\alpha),
\end{equation*}
correspond to isomorphism classes of semisimple representations of
$\Pi_{\lambda}$ of dimension $\alpha$. If $M$ is a semisimple
$\Pi_{\lambda}$-module then we can decompose it into its simple
components $M = M^{\oplus k_1}_1 \oplus \dots \oplus M^{\oplus
k_r}_r$ where the $M_t$ are non-isomorphic simples. If $\beta^{(t)}$
is the dimension vector of $M_t$, then we say $M$ has representation
type \[ (k_1,\beta^{(1)}; \ldots ; k_r,\beta^{(r)}).
\]
One can therefore stratify $\mathcal{N}(\lambda, \alpha)$ by
representation type.

\subsection{}\label{intro}\label{Qleaves} The variety $\mathcal{N}(\lambda, \alpha)$ occurs as the quotient of a fibre of the
moment map and therefore has the structure of a Poisson variety. Any
Poisson variety decomposes into the disjoint union of symplectic
leaves, see \cite{brgo}. For $\mathcal{N}(\lambda, \alpha)$ we show
that this has a straightforward description.

\begin{trm}
The stratification of $\mathcal{N}(\lambda, \alpha)$ by
representation type equals the stratification by symplectic leaves.
\end{trm}
We note that one can describe $\mathcal{N}(\lambda, \alpha)$ as a
version of the quiver varieties introduced by Nakajima, Remark
\ref{lambda}. Our approach of considering the Poisson geometry of
these reductions has been put into a more general framework through
recent work of \cite{cbeg} and \cite{vdb}.

\subsection{}\label{intro2} We apply Theorem \ref{intro} to symplectic
reflection algebras, these are defined in Section \ref{srasection}.
Let $\Gamma$ be a finite subgroup of $SL(2,\C)$ and let $n>1$ be an
integer. Let $H_{\mathbf{c}}$ be the symplectic reflection algebra
(at parameter $t=0$) for the wreath product $\Gamma_n = S_n \ltimes
\Gamma^n$, and let $Z_{\mathbf{c}}$ denote its centre. The variety
$\mx Z_{\mathbf{c}}$ is Poisson. The connection to quivers is via
the McKay correspondence: to the group $\Gamma$ we associate a
quiver, $Q_{\Gamma}$, by choosing an orientation of its McKay graph.
The quiver, $Q$, is obtained from $Q_{\Gamma}$ by adding one vertex
and one arrow from this vertex to an extending vertex of
$Q_{\Gamma}$. We establish the following.

\begin{trm}
For the quiver $Q$, a dimension vector $\alpha$ and a parameter
$\lambda \in \C^I$ depending on $\mathbf{c}$
\[ \mathcal{N}(\lambda, \alpha) \cong \mx Z_{\mathbf{c}}.
\]
This isomorphism is Poisson up to nonzero scalar multiple and in
particular maps symplectic leaves to symplectic leaves.
\end{trm}

Crawley-Boevey and Holland proved this theorem for finite subgroups
of $\mathrm{SL}(2,\C)$, see \cite[Theorem 0.2]{cbh}, although they
do not mention Poisson structures, and our result can be seen as an
extension of theirs. In fact, our proof is based on a result of
Etingof and Ginzburg \cite[Theorem 11.16]{eg} who prove Theorem
\ref{sra-cmiso} with the extra condition that the parameter
$\mathbf{c}$ is generic. The definition of the dimension vector
$\alpha$ and deformation parameter $\lambda$ used in the above
theorem are given in \ref{mckaygraph} where they are written
$\ep_{\infty} + n\delta$ and $\lambda'(\mathbf{c})$ respectively.

\subsection{} The proof of Theorem \ref{intro2} is rather complicated
and involves several stages. The Calogero-Moser space,
$X_\mathbf{c}//\gp(n\delta)$, defined in \ref{cmspace}, plays an
important role as does its equivalent description
$\hat{X}_\mathbf{c}//\ggp(n\delta)$, see \ref{hatX_c}. We can break
down the isomorphism of Theorem \ref{intro2} into the sequence of
isomorphisms. In the diagram below the superscripts indicate the
subsection in which the corresponding isomorphism is proved:
\[\xymatrix{ \mathcal{N}(\lambda, \alpha) \ar[r]^{\! \! \! \! \! \! \! \ref{lincalmwrediso}}_{\! \! \! \! \! \! \! \sim} &
\hat{X}_\mathbf{c}//\ggp(n\delta) \ar[r]^{\! \ref{cmiso}}_{\! \sim}
& X_\mathbf{c}//\gp(n\delta) \ar[r]^{\ \ \ref{sra-cmiso}}_{\ \ \sim}
& \mx Z_\mathbf{c}. }\] All of the above maps are Poisson except the
final one which is only Poisson up to nonzero scalar multiple.

\subsection{} The symplectic leaves of $\mx Z_\mathbf{c}$ are an important
in the representation theory of $H_{\mathbf{c}}$. Brown and Gordon
proved in \cite[Theorem 4.2]{brgo} that if two maximal ideals
$\mathfrak{m}$ and $\mathfrak{n}$ of $Z_{\mathbf{c}}$ lie in the
same symplectic leaf then the factor algebras
$H_{\mathbf{c}}/\mathfrak{m}H_{\mathbf{c}}$ and
$H_{\mathbf{c}}/\mathfrak{n}H_{\mathbf{c}}$ are isomorphic.
Therefore understanding the leaves is useful in studying the
irreducible representations of $H_{\mathbf{c}}$. We have the
following consequence of Theorems \ref{intro} and \ref{intro2}

\begin{thm} There is a one-to-one correspondence between the
symplectic leaves of $\mx Z_{\mathbf{c}}$ and the representation
types of $\mathcal{N}(\lambda, \alpha)$.
\end{thm}

Therefore, using Crawley-Boevey's work on root systems and dimension
vectors of simple modules for $\Pi_{\lambda}$, \cite{cb1}, one can
use linear algebra to calculate information about the symplectic
leaves of $\mx Z_{\mathbf{c}}$, see \ref{leafcalc}.

\subsection{}We organise our paper as follows. In Section
\ref{forms} we recall facts about Poisson varieties and symplectic
leaves. The moment map and the Marsden-Weinstein reduction procedure
is discussed in Section \ref{momentsection}. We examine the special
case of reductions for symplectic vector spaces with a
hyper-K\"{a}hler structure in Section \ref{strats}, identifying the
symplectic leaves in Proposition \ref{leavesforlocal}. We then apply
these to results to representations of deformed preprojective
algebras in Section \ref{RepQ} to prove Theorem \ref{intro}. The
proof of Theorem \ref{intro2} occupies Sections \ref{cmsection} and
\ref{srasection} in which we also introduce Calogero-Moser space. We
calculate the example of the centre of the rational Cherednik
algebra of type $B_n$ in Section \ref{examples}, computing the
parameters at which $\mx Z_{\mathbf{c}}$ is singular and also the
number of symplectic leaves and their dimensions.

\section{Poisson varieties and symplectic leaves}\label{forms}

\subsection{Poisson algebras}\label{poissonvarieties}
Throughout we work over the complex numbers. Thus varieties and
manifolds are complex, and symplectic form will mean complex
symplectic form, unless stated otherwise. In this section we
introduce the notions of Poisson varieties and manifolds and the
stratification by symplectic leaves.

Let $A$ be a commutative $\C$-algebra. $A$ is a \textit{Poisson
algebra} if there exists a $\C$-bilinear bracket $\{-,-\}: A\times A
\to A$ such that

\begin{enumerate}
\item $(A, \{-,-\})$ is a Lie algebra;

\item $\{-,-\}$ satisfies the product rule, that is, $\{a,bc\} =
\{a,b\}c + b\{a,c\}$ for all $a,b,c \in A$.
\end{enumerate}

Let $I$ be an ideal of a Poisson algebra $R$. Then $I$ is a
\textit{Poisson ideal} if $\{R,I\} \subseteq I$. A homomorphism
$\phi: A \to B$ between Poisson algebras $(A, \{-,-\}_A)$ and
$(B,\{-,-\}_B)$ is called \textit{Poisson} if $\phi\{ a_1,a_2\}_A =
\{ \phi(a_1),\phi(a_2)\}_B$ for all $a_1,a_2 \in A$.

\subsection{Poisson varieties} We say that an algebraic variety, $X$, with structure
sheaf $\Co$, is a \textit{Poisson variety} if $\mathcal{O}$ is a
sheaf of Poisson algebras such that the restriction homomorphisms
are Poisson. As an example, given any finitely generated Poisson
algebra $A$, let $\mx A$ denote its set of maximal ideals endowed
with the structure of an algebraic variety. Then $\mx A$ is a
Poisson variety, as is any open subvariety of $\mx A$. We say that a
morphism $\psi : X \to Y$ between Poisson varieties $X$ and $Y$ is
\textit{Poisson} if for all open $U \subseteq Y$, $\psi^{\sharp}_U :
\Co_Y (U) \to \Co_X (\psi^{-1}(U))$ is Poisson. If $X$ and $Y$ are
both affine varieties then $\psi$ is Poisson if and only if the
comorphism $\psi^{\sharp} : \Co_Y (Y) \to \Co_X (X)$ is a Poisson
homomorphism. For any $f\in \Co (X)$ the \textit{Hamiltonian vector
field of f} is the derivation, $\{f,-\}$, of $\Co (X)$.

\subsection{}\label{sympvariety}Let $M$ be a smooth affine algebraic variety, and denote its structure sheaf
by $\mathcal{O}$. We say that $M$ is a \textit{symplectic variety}
if it comes equipped with nondegenerate closed $2$-form, $\omega$. A
symplectic variety $(M,\omega)$ gives rise to a Poisson bracket on
$\mathcal{O} (M)$ as follows. Let $\vect (M)$ denote the Lie algebra
of vector fields on $M$. Nondegeneracy of $\omega$ means that there
is a $\C$-linear map $\mathcal{O} (M) \to \vect (M); f \mapsto
\Xi_f$ where $\Xi_f$ is the unique vector field such that $\omega
(\Xi_f,-) = df$. For all $f,g \in \mathcal{O} (M)$ we define the
bracket of $f$ and $g$ by $\{f,g\} = \Xi_f g = \omega (\Xi_f,
\Xi_g)$. This is a Poisson bracket by \cite[Proposition
5.5.3]{marsratiu}. Hence $M$ is a Poisson variety and $\Xi_f$ is the
Hamiltonian vector field of $f$, see \ref{poissonvarieties} above.

\subsection{Poisson manifolds}\label{poissonmanifolds} For
a manifold, $M$, we denote the algebra of holomorphic functions
defined on $M$ by $\mathcal{C} (M)$. We say that $M$ is
\textit{Poisson manifold} if $\mathcal{C} (M)$ is a Poisson algebra.
A smooth map $f: M \to N$ between Poisson manifolds is called
\textit{Poisson} if the corresponding algebra homomorphism $f^* :
\mathcal{C} (N) \to \mathcal{C} (M)$ is a Poisson homomorphism. Any
symplectic manifold is a Poisson manifold - the argument for showing
that symplectic varieties are Poisson varieties works for manifolds
also. For any $f\in \mathcal{C}(M)$ we define the
\textit{Hamiltonian vector field} of $f$ to be the derivation of
$\mathcal{C}(M)$ given by $\{f,-\}$.

A smooth Poisson variety, $X$, is naturally a manifold. The
existence of a Poisson bracket on $\Co (X)$ is equivalent to the
existence of a bivector $\theta: X \to T^2X$ which is a section to
the bundle map and satisfies $[\theta, \theta] = 0$. Here $TX$ is
the tangent bundle of $X$ and $[-,-]$ is the Schouten bracket.
Considering the bivector as complex analytic yields a Poisson
bracket on $\mathcal{C} (X)$, making $X$ a Poisson manifold.

\subsection{}\label{leaffacts}
Let $M$ be a Poisson manifold.
\begin{defn}
The \textit{symplectic leaf} $\mathcal{S}(p)$ containing a point $p$
of $M$ is the set of points $q$ which are connected to $p$ by
piecewise smooth curves, each segment of which is the integral curve
of a Hamiltonian vector field.
\end{defn}

Therefore to find the symplectic leaf containing $p$ one works out
all the points one can reach by travelling along integral curves of
Hamiltonian vector fields at $p$. For each such point, $q$,
connected to $p$ in this way one repeats the process, finding the
points one can reach by travelling along integral curves of
Hamiltonian vector fields at $q$. Then one continues this process
until the symplectic leaf is swept out.

\subsection{} The \textit{rank} of $p$ is equal to the dimension of the
subspace of $T_{p} M$ spanned by the Hamiltonian vector fields
evaluated at $p$. For Poisson manifolds $M$ and $N$ such that $N
\subseteq M$ we shall say that $N$ is a \textit{Poisson submanifold}
of $M$ if the inclusion map is Poisson.

The proof of the following theorem, under some additional
hypotheses, goes back to Lie, and the general case is due to
Kirillov, \cite{Kirillov76}; for a proof in the complex case see
\cite[Theorem 3.26]{Adlerbook}. Further discussion about symplectic
leaves can be found in \cite[Section 10]{marsratiu} and in \cite{w}.

\begin{thm}
$M$ is a disjoint union of its symplectic leaves. Each leaf is a
symplectic manifold which is a Poisson submanifold of $X$ and the
dimension of the leaf through $p$ is equal to the rank of $p$ in
$M$.
\end{thm}

\subsection{}Suppose that $X$ is a smooth affine Poisson variety. Then,
as described in \ref{poissonmanifolds}, $X$ is a Poisson manifold so
one can stratify $X$ by symplectic leaves. Although $X$ is an
algebraic variety it is quite possible that the leaves are not
algebraic varieties, see \cite[Remarks 3.6 (1)]{brgo}.

Now suppose that $X$ is not necessarily smooth; we can stratify $X$
by symplectic leaves, as described in \cite[$\S 3.5$]{brgo}, as
follows. Let $U_0$ be the smooth locus of $X$. Then, since $X$ is
affine, $U_0$ is a smooth Poisson variety and so is a Poisson
manifold. We can stratify $U_0$ by symplectic leaves, say $U_0 =
\bigsqcup_{i\in \mathcal{I}_0} \mathcal{S}_{i,0}$. Now we proceed
inductively by setting $X_0 = X$ and defining $X_k = X_{k-1}
\setminus U_{k-1}$ for $k\geq 1$. $X_k$ is an affine Poisson variety
by \cite[Proposition 15.2.14(i)]{mcrob}, so one can, as above,
stratify the smooth locus, $U_k$, of $X_k$ by symplectic leaves,
$U_k = \bigsqcup_{i\in \mathcal{I}_k} \mathcal{S}_{i,k}$. Then
\begin{equation}X = U_0 \sqcup \dots \sqcup U_t\end{equation} for
some $t$ and we call
\begin{equation*}
X = \bigsqcup_{\substack{i\in \mathcal{I}_k,\\0\leq k \leq
t}}\mathcal{S}_{i,k}
\end{equation*}
the \textit{stratification of X by symplectic leaves}.

\subsection{}\label{algleaves}We say that a closed subvariety, $Y$, of $X$ is
\textit{Poisson} if the inclusion $Y \hookrightarrow X$ is Poisson.
This is equivalent to the condition that the defining ideal of $Y$
is a Poisson ideal.

\begin{prop}\cite[Proposition 3.7]{brgo}
If $Y$ is a closed Poisson subvariety of $X$ and $\s$ is a
symplectic leaf of $X$ then $\s \cap Y \neq \emptyset$ implies that
$\s \subseteq Y$. Furthermore, if $X$ is a finite union of
symplectic leaves then each leaf is an irreducible locally closed
smooth subvariety of $X$.
\end{prop}

\section{The moment map}\label{momentsection}
\subsection{}\label{moment}
For details of the following see \cite[$\S 1.4$]{chrissgin}. Let
$(M,\omega)$ be a symplectic variety as in \ref{sympvariety} and
suppose that a reductive algebraic group $G$ acts by morphisms on
$M$ preserving the symplectic form. Let $\mathfrak{g}$ be the Lie
algebra of $G$.

A vector field $X$ is \textit{symplectic} if it preserves the
symplectic form on $M$, that is, $L_X \omega = 0 $, where $L_X$ is
the Lie derivative with respect to $X$. We denote by $\symp(M)$ the
Lie subalgebra of symplectic vector fields on $M$.

\begin{prop}\cite[Proposition 1.2.5]{chrissgin}
$\Xi_f$ is a symplectic vector field for all $f \in \mathcal{O}(M)$,
and the map $f \mapsto \Xi_f$ defines a Lie algebra homomorphism
$(\mathcal{O}(M), \{ -,-\} ) \to (Symp(M), [-,-] )$.
\end{prop}

Recall the infinitesimal action of $\mathfrak{g}$ on
$\mathcal{O}(M)$: for each $x\in \mathfrak{g}$ and $f\in \Co (M)$ ,
$x_M (f) := \frac{d}{dt} (exp(tx)\circ f) |_{t=0}$. The operator
$x_M$ is a derivation of $\Co (M)$ and so defines a map
$\mathfrak{g} \to \mathrm{Vect}(M)$. This map is a Lie algebra
homomorphism, \cite[Proposition 9.3.6]{marsratiu}, and one can
easily check that its image is contained in the space of symplectic
vector fields.

The $G$-action is said to be \textit{Hamiltonian} if there exists a
Lie algebra homomorphism $H: \mathfrak{g} \to (\mathcal{O} (M),
\{-,-\})$ such that the following diagram commutes:

$$\xymatrix{ \mathfrak{g} \ar[dr] \ar_H[d] &      \\
\mathcal{O}(M) \ar[r] &      \symp(M).}$$

Suppose that the $G$-action is Hamiltonian and let $H_x = H(x)$ for
all $x\in \mathfrak{g}$. The \textit{moment map} $\mu : M \to
\mathfrak{g}^{\ast}$ is the morphism of algebraic varieties defined
by $\mu (m) (x) = H_x (m)$ for all $m\in M$ and $x\in \mathfrak{g}$.

\begin{example}
\begin{rm}Recall that the coordinate ring, $\Co (\g^{\ast})$, of
$\g^*$, equals the symmetric algebra  $S \mathfrak{g}$ of
$\mathfrak{g}$. The algebra $S \g$ has a canonical Poisson bracket
induced from the Lie bracket on $\mathfrak{g}$. One can calculate
the infinitesimal action of $\g$ on $S\g$ (see \cite{marsratiu},
Examples (c), Page 272): for $x \in \g$ and $f \in \g$,
$x_{\g^{\ast}}(f) = [x,f]$. The action of $G$ is Hamiltonian, with
$H_x = x$ for all $x\in \g$, and the corresponding moment map is
$\mu = \mathrm{Id}_{\g^{\ast}}$.

Let $\el$ be a closed orbit (under the coadjoint action) in
$\g^{\ast}$. Let $I$ be the defining ideal of $\el$. Then $I$ is a
Poisson ideal of $S \g$ so that $\el$ is a Poisson variety. In fact
this makes $\el$ a symplectic variety, \cite[Proposition
1.3.21]{chrissgin}. Now for all $x\in \g$, $x_{\el}$ is the
derivation of $\Co (\el) = S\g/I$ given by $x_{\el}(f + I) = [x,f] +
I$ for all $f\in \g \subset S\g$. One sees that the action of $G$ on
$\el$ is Hamiltonian with $H_x = x + I$ and the corresponding moment
map is simply the inclusion $\el \hookrightarrow \g^*$.
\end{rm}
\end{example}

\subsection{}The moment map has the following properties.

\begin{prop}\cite[Lemma 1.4.2]{chrissgin}
The map
\[ \mu^{\sharp} : \C [\mathfrak{g}] \to \Co (M) \]
induced by $\mu$ is a Poisson homomorphism, and if $G$ is connected
then $\mu$ is $G$-equivariant (relative to the coadjoint action on
$\mathfrak{g}^{\ast}$).
\end{prop}

\subsection{}In general it is interesting to know whether a group
action on a symplectic variety is Hamiltonian. In the linear case
this is known.

\begin{trm}\cite[Proposition 1.4.6]{chrissgin}\label{mw} Let $V$ be
a symplectic vector space with symplectic form $\omega$ and let $G$
be a reductive algebraic subgroup of $\mathrm{Sp} V$. Then the
action is Hamiltonian with $H_A(v) = \frac{1}{2} \omega(A\cdot v,
v)$ for all $A\in \mathrm{Lie}\ G$ and $v\in V$ and the
corresponding moment map is $G$-equivariant.
\end{trm}

\subsection{Reduction}\label{poisson}
Many of the varieties which appear in this paper arise as quotient
spaces of fibres of a moment map, and we shall see now how these
spaces carry a Poisson structure. Let $(M,\omega)$ be a symplectic
variety and let the reductive algebraic group $G$ act by morphisms
on $M$ preserving $\omega$. Suppose this action is Hamiltonian and
that the corresponding moment map $\mu$ is $G$-equivariant. Let
$\mathcal{L}$ be a closed orbit under the coadjoint action of $G$ on
$\mathfrak{g}^{\ast}$, with defining $G$-stable ideal $J \lhd
\C[\mathfrak{g}]$. Then $\mu^{-1}(\mathcal{L})$ is a $G$-stable
closed subvariety of $M$, with defining $G$-stable radical ideal $I
\lhd \Co (M)$.

We consider the quotient variety $M_{\el} := \mu^{-1}(\mathcal{L})
// G$, which we call a \textit{Marsden-Weinstein
reduction}. The double lines indicate that we consider closed
orbits, in other words points in the maximal ideal spectrum of $\Co
(\mu^{-1}(\el))^G$. Let $\{-,-\}$ be the bracket on $\Co (M)$. Since
$G$ is reductive $\big( \Co (M)/I \big)^G \cong \Co (M)^G/ I^G$ so
one can define a bracket on $\big( \Co (M)/I \big)^G$ by defining
one on $\Co (M)^G/ I^G$. Define a bracket, $\{-,-\}'$, on $\Co
(M)^G/ I^G$ by $\{f + I^G, g + I^G\}' = \{f,g\} + I^G$ for all
$f,g\in \Co (M)^G$.

\begin{prop}
The bracket $\{-,-\}'$ is well-defined and is a Poisson bracket on
$\Co (M_{\el})$.
\end{prop}
\begin{proof}
The bracket $\{-,-\}'$ will clearly yield a Poisson bracket as long
as it is well defined. To see this we note first that $\{f,g\} \in
\Co (M) ^G$ for all $f,g \in \Co (M)^G$. It remains to show that
$I^G$ is a Poisson ideal of $\Co (M)^G$.

For all $x\in \mathfrak{g}$, $\mu^{\sharp}(x) = x\circ \mu$ and if
we evaluate this at an element of $m\in M$ we see that $(x\circ
\mu)(m) = x (H_{-}(m)) = H_x (m)$. In short, $\mu^{\sharp}(x) = H_x$
and it follows that $I' := \Co(M)\mu^{\sharp}(J)$ is generated by
polynomials in the $H_x$.

Now, if $f\in \Co (M)^G$ then $h\cdot f = f$ for all $h\in G$ and
therefore $x_M f = 0$ for all $x\in \mathfrak{g}$. Thus $\{H_x,f\} =
x_M f = 0$ and, by the product rule, $\{f,i\} \in I'$ for all $i\in
I'$. Therefore $I'^G$ is a Poisson ideal of $\Co (M)^G$. Finally, $I
= \sqrt{I'}$ implies that $I^G = \sqrt{I'^G}$, and so $I^G$ is a
Poisson ideal of $\Co (M)^G$ by \cite[3.3.2]{dix}.
\end{proof}

\begin{example}[\textbf{Rank one matrices}]\begin{rm}
Let $M_n = \mat(n\times n, \C)$ which is the Lie algebra of
$\gl(n,\C)$. $M_n$ can be identified with with its dual via the
trace pairing and in this way becomes a Poisson variety as in
Example \ref{moment}. One can describe $\Co(M_n)$ as $\C[x_{ij}: 1
\leq i,j \leq n]$ and then its Poisson bracket is $\{x_{ij} ,
x_{kl}\} = - \delta_{jk} x_{il} + \delta_{il}x_{kj}$.

Let $R_l$ be the subvariety of $M_n$ consisting of matrices whose
rank is less than or equal to $l$. Thus $R_l$ is the set of matrices
whose $(l+1) \times (l+1)$ minors are zero and so is a closed
subvariety.
The defining ideal of $R_l$ is the prime ideal, $I_l$, of
$\C[x_{ij}]$ generated by the generic $(l+1)\times (l+1)$ minors,
\cite[Theorem 2.10]{Brunsbook}. We are interested in looking at the
matrices whose rank is less than or equal to one, that is, the
subvariety $R_1$. By \cite[Proposition 15.2.14]{mcrob}, $I_1$ is a
Poisson ideal.

Let $0 \neq k \in \C$ and let $T_k \subset M_n$ be the closed
subvariety of matrices whose trace equals $k$. Then $T_k$ has
defining ideal $(x_{11} +\dots + x_{nn} - k) \lhd \C[x_{ij}]$ and a
direct calculation shows that this a Poisson ideal. From this it
follows that $U_k := R_1 \cap T_k$ is a closed Poisson subvariety of
$M_n$ with defining ideal $J_k := \sqrt{I_1 + (x_{11} +\dots +
x_{nn} - k)}$.

We can describe $U_k$ as a reduction for a moment map, as follows.
Let $V$ be an $n$-dimensional complex vector space and consider the
symplectic space $W = V \oplus V^*$. Using the usual scalar product
on $V$ we can think of $V$ as consisting of column vectors and $V^*$
of row vectors. We fix bases $c_1, \dots ,c_n \in V$ and $r_1 ,
\dots , r_n \in V^*$. There is an action of $\C^{\times}$ by
$\lambda \cdot (c,r) = (\lambda v , \lambda^{-1} r)$ for all
$\lambda \in \C^{\times}, c\in V$ and $r \in V^*$. This action is
Hamiltonian, Theorem \ref{mw}, with moment map $\mu : W \to \C;\
\mu(c,r) = rc$, here we identify $\C$ with $\C^* = \lie
(\C^{\times})^*$ via the trace pairing. It is an easy calculation to
see that \[\Co (W)^{\C^{\times}} = (\Co (V) \otimes \Co
(V^*))^{\C^{\times}} = \C [ r_i\otimes c_j: 1 \leq i,j \leq n]\] and
then it follows that \[\Co (\mu^{-1}(k)//\C^{\times}) = \C [
r_i\otimes c_j]/I_k\] where $I_k := \sqrt{(r_1\otimes c_1 + \dots +
r_n\otimes c_n - k)} \lhd \C [ r_i\otimes c_j]$. One can check
directly that $I_k$ is a Poisson ideal of $\C [ r_i\otimes c_j]$,
which verifies Proposition \ref{poisson} in this example.

We now describe the isomorphism between $\mu^{-1}(k)//\C^{\times}$
and $U_k$. Let $G = \gl (n, \C)$ which acts naturally on $V$ and by
conjugation on $M_n \cong \lie (G)^*$. One obtains an action of $G$
on $W$ and the corresponding moment map (Theorem \ref{mw}) is \[ \nu
: W \to M_n;\ (r,c) \mapsto r \otimes c \in V \otimes V^* \cong
M_n.\] It is clear that the image of $\nu$ is contained in $R_1$.
The map $\nu$ is equivariant for the action of $G$ and therefore
also for the action of $\C^{\times}$, where $\C^{\times}$ is thought
of as the subgroup of nonzero scalar matrices in $G$. Therefore,
since $\C^{\times}$ acts trivially on $M_n$ we get a morphism
induced from $\nu$: \[ t: \mu^{-1}(k)//\C^{\times} \to U_k,\] whose
comorphism is the map \[ \C[x_{ij}]/J_k \to \C[r_i \otimes
c_j]/I_k;\ x_{ij} + J_k \mapsto r_i \otimes c_j + I_k.\] This is a
Poisson isomorphism and so we have described $U_k$ as the reduction
$\mu^{-1}(k)//\C^{\times}$.

\end{rm}
\end{example}

\subsection{}\label{shiftingtrick}Given a reduction over a closed coadjoint orbit, one can
perform the so-called shifting trick to express this as reduction
over a fixed point.

\begin{lem}[\textbf{The shifting trick}] Let $\el$ be a closed
coadjoint orbit of $\g^*$ and let $\xi \in \g^*$ be a fixed point.
Then $M \times -\el$ is a symplectic variety such that the natural
action of $G$ is Hamiltonian with corresponding moment map $\mu': M
\times -\el \to \g^*;\ (m,x) \mapsto \mu(m) + x$. There is an
isomorphism of Poisson varieties $\mu'^{-1}(\xi)//G \cong
\mu^{-1}(\el + \xi)//G$.
\end{lem}
\begin{proof}
As explained in Example \ref{moment}, $-\el$ is a symplectic variety
with Hamiltonian $G$-action such that the moment map is just the
inclusion $-\el \hookrightarrow \g^*$. Therefore the product $M
\times -\el$ is symplectic and the action of $G$ on this product is
Hamiltonian because it is Hamiltonian on each component. To see this
one takes the Hamiltonian maps $H_1, H_2$ of $M$ and $-\el$
respectively and checks that the map $H = H_1 \otimes 1 + 1 \otimes
H_2 : \g \to \Co (M\times -\el) = \Co (M) \otimes \Co (-\el)$ is a
Hamiltonian for the action of $G$ on the product. The moment map
corresponding to $H$ is then the sum of the moment maps for $M$ and
$-\el$. The projection map $M \times -\el \to M$ is a
$G$-equivariant Poisson map and restricts to an isomorphism between
$\mu'^{-1}(\xi)$ and $\mu^{-1}(\el + \xi)$. The $G$-equivariance of
the projection map means that this induces an isomorphism
$\mu'^{-1}(\xi)//G \cong \mu^{-1}(\el + \xi)//G$. This isomorphism
is Poisson because the projection map is Poisson and by the
proposition above.
\end{proof}

\section{The local normal form for the moment
map}\label{strats}

\subsection{}\label{hyperkahlerdefn} We discuss a very particular
case of reduction which will apply to representations of
preprojective algebras. Our arguments follow closely those of
\cite{sjaler} and \cite[Section 3]{nak1}. Let $M$ be a finite
dimensional complex vector space. Let $G$ be a complex connected
reductive algebraic group acting linearly on $M$, and we think of
$G$ as the complexification of a compact real algebraic group $R$.
Let $I : M \to M$ be the map defined by multiplication by
$\sqrt{-1}$ and let $(-,-)$ be a real inner product on $M$ which is
invariant under the actions of $I$ and $R$. Thus $M$ is a K\"{a}hler
manifold. Let $J,K: M \to M$ be $\R$-linear maps preserving $(-,-)$
and satisfying $I^2 = J^2 = K^2 = IJK = -\mathrm{Id}$. This makes
$M$ into a hyper-K\"{a}hler manifold, \cite{hitchin}.

\subsection{}\label{hyperkahlerforms} We have real symplectic forms on
$M$ given by $\omega_1 (v,w) = (v,Iw)$, $\omega_2 (v,w) = (v,Jw)$
and $\omega_3 (v,w) = (v,Kw)$ for all $v,w \in M$. If we define
$\omega = \omega_2 + I\omega_3$ then this is a complex symplectic
form on $M$. We shall assume that the action of $G$ preserves
$\omega$. Let $\mathfrak{r}, \mathfrak{g}$ be the Lie algebras of
$R$ and $G$ respectively. We have two moment maps to consider in
this situation. The first, $\mu_1 : M \to \mathfrak{r}^*$, is a
moment map with respect to $\omega_1$ and is obtained completely
analogously to that of Theorem \ref{mw}. We also have the moment map
$\mu : M \to \mathfrak{g}^*$ with respect to $\omega$ which is
constructed as in Theorem \ref{mw}.

\subsection{}\label{localnorm} We wish to consider a closed
orbit corresponding to a point in $\mu^{-1}(0)//G$, say $Gm$, for
some $m\in \mu^{-1}(0)$. Let $H$ be the stabiliser of $m$ in $G$,
and denote its Lie algebra by $\h$. Let $T_m (Gm)$ be the tangent
space of $m$ in $Gm$, and let $(T_m (Gm))^{\omega} = \{v\in T_m M\
:\ \omega(v,w)=0\ \mathrm{for\ all}\ w\in T_m (Gm)\}$. By \cite[page
324]{gust} $T_m (Gm)\subseteq (T_m (Gm))^{\omega}$; set $\hat{M} =
(T_m (Gm))^{\omega}/T_m (Gm)$. The action of $H$ preserves the
induced symplectic structure on $\hat{M}$, let $\hat{\mu}: \hat{M}
\to \h^*$ be the corresponding $H$-equivariant moment map, Theorem
\ref{mw}.

We choose an $\mathrm{Ad}(H)$-invariant splitting $\g = \h \oplus
\h^{\bot}$ and a dual splitting $\g^* = \h^* \oplus \h^{\bot *}$.
There is a natural action of the reductive group, $H$, on
$\mathfrak{h}^{\bot *} \times \hat{M}$ by $h\cdot (\xi, \hat{m}) =
(\mathrm{Ad}^*(h)\xi, h\hat{m})$. The associated bundle $(G \times
\h^{\bot *} \times \hat{M})//H$ is a symplectic variety with
Hamiltonian $G$-action (see \cite[Page 164]{nak1}) - the $G$-action
is given by $x\cdot [g,\xi , \hat{m}] = [xg, \xi, \hat{m}]$ for all
$x\in G$. The ($G$-invariant) moment map is given by
\[ J(H\cdot(g, \xi, \hat{m})) = \mathrm{Ad}^* (g)(\xi +
\hat{\mu}(\hat{m})). \]

\subsection{}\label{localnormthm}We shall refer to the following theorem by saying that
there exists a \textit{local normal form for} $\mu$. Let $\Upsilon =
(G \times \h^{\bot *} \times \hat{M})//H$.
\begin{thm}\cite[Lemma 3.2.1]{nak1}
Let $\pi$ and $\rho$ be the $G$-orbit maps for $M$ and $\Upsilon$
respectively. Let $y = \pi(m)$. Then there exist complex
neighbourhoods, $U$, of $y$ and $W$ in $(\h^{\bot *} \times
\hat{M})//H$ and biholomorphisms $\Psi, \tilde{\Psi}$ such that the
following diagram commutes \[ \xymatrix{ \pi^{-1}(U) \ar_{\pi}[d]
\ar^{\tilde{\Psi}}[r] & \rho^{-1}(W) \ar^{\rho}[d] \\ U
\ar_{\Psi}[r] & W.}\] Furthermore, $\tilde{\Psi}$ is a
$G$-equivariant map which intertwines symplectic forms and maps $m$
to $(H, 0, 0) \in \Upsilon$.
\end{thm}



\subsection{Stratifying by orbit type}\label{orbittype}

Let $\el$ be a closed coadjoint orbit in $\mathfrak{g}^{\ast}$ and
let $Z = \mu^{-1}(\el)$. Let $M_{\el}$ denote the corresponding
Marsden-Weinstein reduction, $Z//G$. Let $\pi : Z \to M_{\el}$
denote the orbit map.

One can stratify $M_{\el}$ in a natural way. For $z \in Z$ let $G_z$
be the stabiliser of $z$ in $G$. Let $\mathcal{T}$ be the set of
conjugacy classes of subgroups $G_z$ in $G$. For $\phi \in M_{\el}$
denote by $T(\phi)$ the conjugacy class of stabilisers belonging to
the unique closed orbit in $\pi^{-1}(\phi)$. Define
$(M_{\el})_{\tau} = \{\phi \in M_{\el} : T(\phi) = \tau\}$. We call
$M_{\el} = \bigcup_{\tau \in \mathcal{T}}(M_{\el})_{\tau}$ the
\textit{stratification by orbit type}. For each $\tau \in
\mathcal{T}$ choose a stabiliser, $G_{\tau}$, representing the
conjugacy class $\tau$. For all $\tau_1, \tau_2 \in \mathcal{T}$ we
write $G_{\tau_1} \leq_c G_{\tau_2}$ if $G_{\tau_1}$ is conjugate to
a subgroup of $G_{\tau_2}$. We define a partial order on
$\mathcal{T}$ by $\tau_1 \geq \tau_2$ if $G_{\tau_1} \leq_c
G_{\tau_2}$.

The proposition below is due to Schwarz, another proof can be found
in \cite{domzub}.

\begin{prop}\cite[Lemma 5.5]{schwarz}\label{orbittypelocclosed}\mbox{}
For each $\tau \in \mathcal{T}$ the stratum $(M_{\el})_{\tau}$ is
irreducible and locally closed and its closure is $\bigcup_{\nu \leq
\tau}(M_{\el})_{\nu}$. The stratification $M_{\el} = \bigcup_{\tau
\in \mathcal{T}} (M_{\el})_{\tau}$ is finite.
\end{prop}

We would like to compare the stratification by orbit type to the
stratification by symplectic leaves. In particular that would
require that each stratum is a symplectic manifold.

\subsection{}We shall need a useful lemma. For any $H \leq G$ let $Z_{(H)} = \{
z\in Z : G_z\ \mathrm{is\ conjugate\ to}\ H\}$. For any $\tau \in
\mathcal{T}$ such that $(M_{\el})_{\tau}$ is nonempty, we choose
some $H\in \tau$ and consider $Z_{(H)}$. The restriction of $\pi$ to
$Z_{(H)}$ is a map whose image contains $(M_{\el})_{\tau}$.

\begin{lem}\label{nastylemma}
Let $\tau \in \mathcal{T}$ and let $H \in \tau$. Suppose that
$(M_{\el})_{\tau} \neq \emptyset$. Then
$\pi^{-1}((M_{\el})_{(\tau)}) \cap Z_{(H)}$ is a nonempty open
$G$-stable subset of $Z_{(H)}$.
\end{lem}
\begin{proof}
Let $z\in Z_{(H)}$. There exists $\sigma \in \mathcal{T}$ such that
$\pi (z) \in (M_{\el})_{\sigma}$. Let $X_{\sigma} = \bigcup_{\nu
\leq \sigma}(M_{\el})_{\nu}$ which is a closed subset of $M_{\el}$
by the proposition above. If the orbit $G \cdot z$ is closed in $Z$
then $\sigma = \tau$ by definition of $(M_{\el})_{\sigma}$. Suppose
that this orbit is not closed. Then $\sigma \neq \tau$ and the
unique closed orbit in $\pi^{-1}(\pi(z))$ has dimension strictly
less than $\mathrm{Dim} G\cdot z$ by \cite[Corollary
13.3.1]{Santosbook}. Let $K\in \sigma$. By comparing dimensions we
see that $K \nleq_c H$ and this implies that $\tau \nleq \sigma$.
Therefore $(M_{\el})_{\tau} \cap X_{\sigma} = \emptyset$.

Let $\s = \{ \sigma \in \mathcal{T}: \pi(z) \in (M_{\el})_{\sigma}\
\mathrm{for\ some}\ z\in Z_{(H)}\}$ and let $\s' = \s \setminus
\{\tau \}$. By definition we have $\pi^{-1}(\bigcup_{\sigma \in
\s}X_{\sigma}) \cap Z_{(H)} = Z_{(H)}$. By the proposition above
$X_{\tau}\setminus (M_{\el})_{\tau}$ is closed in $M_{\el}$.
Therefore $\bigcup_{\sigma \in \s}X_{\sigma}\setminus
(M_{\el})_{\tau} =  (X_{\tau}\setminus (M_{\el})_{\tau}) \cup
\bigcup_{\sigma \in \s'}X_{\sigma}$ is closed in $M_{\el}$ and
\[\begin{split}\pi^{-1}((M_{\el})_{(\tau)}) \cap Z_{(H)} &= \big(
\pi^{-1}(\bigcup_{\sigma \in \s}X_{\sigma}) \cap Z_{(H)} \big)
\setminus \big( \pi^{-1}(\bigcup_{\sigma \in \s}X_{\sigma}\setminus
(M_{\el})_{\tau}) \cap Z_{(H)}\big)\\ &= Z_{(H)} \setminus \big(
\pi^{-1}(\bigcup_{\sigma \in \s}X_{\sigma}\setminus
(M_{\el})_{\tau}) \cap Z_{(H)}\big)\end{split}\] is open in
$Z_{(H)}$. That this set is nonempty and $G$-stable follows from the
definition of $Z_{(H)}$ and the fact that $\pi$ is constant on
$G$-orbits.
\end{proof}

\subsection{}
We prove a result crucial to the proof of Theorem \ref{intro}. The
following theorem is based very closely on \cite[Theorem
2.1]{sjaler}. This earlier result was proved over $\R$ and for a
compact Lie group, we provide full details to verify that the result
carries over to our situation.

\begin{thm}\label{strat}
Let $\el$ be a fixed point in $\g^*$. Then the decomposition
$M_{\el} = \bigsqcup \{ (M_{\el})_{\tau}:\ {\tau} \in \mathcal{T}\}$
is a stratification of $M_{\el}$ into a disjoint union of symplectic
manifolds. Let $\tau \in \mathcal{T}$ and choose $H\in \tau$. Then
the pullback of the symplectic form, $(\omega_0)_{\tau}$, of
$(M_{\el})_{\tau}$ to $\pi^{-1}((M_{\el})_{\tau}) \cap Z_{(H)}$
equals the restriction to $\pi^{-1}((M_{\el})_{\tau}) \cap Z_{(H)}$
of the symplectic form $\omega$.
\end{thm}
\begin{proof}

We first give the proof when $\el = \mathbf{0}$. Let $\zeta \in
(M_{\el})_{\tau}$ and $G m \subseteq Z$ be the unique closed orbit
in $\pi^{-1}(\zeta)$. Let $\Upsilon = (G\times \h^{\bot *} \times
\hat{M})//H$. By Theorem \ref{localnormthm} we can work in the model
space $\Upsilon$, that is, there is $G$-equivariant biholomorphism,
intertwining symplectic structures, between a neighbourhood of $Gm$
in $M$ and a neighbourhood of $G/H\times 0 \times 0$ in $\Upsilon$.

Let $\Upsilon_{(H)} = \{y\in \Upsilon : G_y\ \mathrm{is\ conjugate\
to}\ H\}$, where $G_y$ denotes the stabiliser of $y$ with respect to
the action of $G$ on $\Upsilon$. Let $y\in \Upsilon_{(H)}$. Then we
can write $y = [g,\xi,\hat{m}]$, and for all $x\in G_y$,
$[xg,\xi,\hat{m}] = [g,\xi,\hat{m}]$, that is, there exists $h\in H$
such that $(xg,\xi,\hat{m})=(gh^{-1},\mathrm{Ad}^*(h)\xi,h\hat{m})$.
By considering the first component of $y$ we deduce that $g^{-1}G_yg
= H$, and therefore for all $h\in H$, $\mathrm{Ad}^*(h)\xi = \xi$
and $h\hat{m} = \hat{m}$. It follows from Theorem \ref{mw} that
$h\hat{m} = \hat{m}$ for all $h\in H$ implies $\hat{\mu}(\hat{m}) =
0$. Therefore \[
\begin{split}J^{-1}(0) \cap \Upsilon_{(H)} &= \{[g,\xi,\hat{m}]:
\mathrm{Ad}^*(h)\xi
 = \xi, h\hat{m}=\hat{m}\ \mathrm{for\ all}\ h\in H\ \mathrm{and}\  \xi + \hat{\mu}(\hat{m}) =
0 \}\\ &=\{[g,\xi,\hat{m}]: \xi= 0\ \mathrm{and}\ h\hat{m}=\hat{m}\
\mathrm{for\ all}\ h\in H\}\\ &= (G\times \hat{M}_H)//H = G/H \times
\hat{M}_H, \end{split} \] where $\hat{M}_H$ is the subspace of fixed
points of $\hat{M}$, which is a symplectic subspace of $\hat{M}$.

We note that the map $\tilde{\Psi}$ from Theorem \ref{localnormthm}
maps points of $Z_{(H)}$ to $J^{-1}(0) \cap \Upsilon_{(H)}$. Thus
the calculation of the previous paragraph and the local normal form
imply that $Z_{(H)}$ is a (not necessarily connected) submanifold of
$M$. By Lemma \ref{nastylemma} for any nonempty open $U \subseteq
(M_{\el})_{\tau}$, $\pi^{-1}(U) \cap Z_{(H)}$ is a nonempty open
$G$-stable subset of $Z_{(H)}$. Therefore by Theorem
\ref{localnormthm} we have open sets $U \subseteq (M_{\el})_{\tau}$
and $W \subseteq \hat{M}_H$ and a commutative diagram
\begin{equation}\label{stratasymplectic}
\xymatrix{ \pi^{-1}(U) \cap Z_{(H)} \ar^{\tilde{\Psi}}[r]
\ar_{\pi}[d] & G/H \times W \ar^{p}[d] \\ U \ar_{\Psi}[r] & W}
\end{equation}
where $\Psi, \tilde{\Psi}$ are biholomorphisms, $\tilde{\Psi}$ is
$G$-equivariant and intertwines the restrictions of the symplectic
forms to $Z_{(H)}$ and $G/H \times \hat{M}_H$, and $p$ is the
projection map.

Now using the map $\Psi$ we can give $(M_{\el})_{\tau}$ the
structure of a symplectic manifold. Note that by Proposition
\ref{orbittypelocclosed} $(M_{\el})_{\tau}$ is irreducible and
therefore connected. The statement about the pullback of
$(\omega_0)_{\tau}$ follows from (\ref{stratasymplectic}) since the
pullback by $p$ of the symplectic form on $W$ equals the restriction
of the symplectic form of $\Upsilon$ to $G/H \times W$.

For a general fixed point $\el$ we use the `shifting trick' of Lemma
\ref{shiftingtrick} to reduce to the case above. $\el$ is a
symplectic variety in a trivial way, so we form the product $M
\times -\el$, which is a symplectic variety with moment map
$\mu'(m,-\el) = \mu(m) + -\el$. The local normal form for $\mu'$
will exist because it does so for $\mu$, and we can identify
$\mu'^{-1}(\mathbf{0})//G$ with $\mu^{-1}(\el)//G$.
\end{proof}


\subsection{}\label{leaves} We can define a bracket, $\{-,-\}_1$, on
$\Co (M_{\el})$ in the following way. Let $f,g \in \Co (M_{\el})$
and $p\in (M_{\el})_{\tau}$ for some $\tau \in \mathcal{T}$. Then
$\{f,g\}_1(p) = \{f\big|_{\tau} ,g\big|_{\tau}\}_{\tau} (p)$ where
$f\big|_{\tau} ,g\big|_{\tau}$ denote the restrictions to
$(M_{\el})_{\tau}$ of $f$ and $g$ respectively and $\{-,-\}_{\tau}$
is the Poisson bracket induced by the symplectic form,
$(\omega_0)_{\tau}$. It is not clear that this is a Poisson bracket
on $\Co (M_{\el})$; we remedy this below. Let $d\pi_z$ denote the
differential of $\pi$ at $z\in Z$. For the 2-form
$(\omega_0)_{\tau}$ on $(M_{\el})_{\tau}$ we denote by
$\pi^{*}(\omega_0)_{\tau}$ the pullback of $(\omega_0)_{\tau}$ by
$\pi$. The following proposition is based on \cite[Proposition
3.1]{sjaler}.

\begin{prop}
Let $\{-,-\}_2$ denote the Poisson bracket on $\Co (M_{\el})$
defined in \ref{poisson}. Then $\{-,-\}_1 = \{-,-\}_2$.
\end{prop}
\begin{proof}
Let $f,g \in \Co (M_{\el})$. Let ${\tau}\in \mathcal{T}$, $H\in
\tau$ and $p_0 \in (M_{\el})_{\tau}$. It suffices to show that
$\{f,g\}_1(p_0) = \{f,g\}_2(p_0)$. Let $G p$ be the unique closed
orbit in $\pi^{-1}(p_0)$ so that $p \in Z_{(H)}$. Let $\tilde{f},
\tilde{g} \in \Co (M)^G$ be such that $\tilde{f} + I^G = f$ and
$\tilde{g} + I^G = g$ where $I$ is the ideal of functions in $\Co
(M)$ vanishing on $Z$. Then $\{f,g\}_2(p_0) = \{\tilde{f},
\tilde{g}\}(p)$.

Let $S = \pi^{-1}((M_{\el})_{\tau}) \cap Z_{(H)}$. The proof of
Theorem \ref{strat} shows that $S$ is a submanifold of $M$. It
follows from \cite[Proposition 10.5.2]{marsratiu} that for any $z\in
Z_{(H)}$ travelling along the integral curve to $\Xi_{\tilde{f}}(z)$
preserves stabiliser type, that is, this integral curve is contained
in $M_{(H)}$. Proposition \ref{poisson} implies that this integral
curve is also contained in $Z$. Therefore $\Xi_{\tilde{f}}(z)$ is
contained in the tangent space to $z$ in $S$.

Let $[\tilde{f}|_S]$ denote the germ of $\tilde{f}$ around $p$ in
$S$, and let $[f|_{\tau}]$ denote the germ of $f$ around $p_0$ in
$(M_{\el})_{\tau}$. By Theorem \ref{strat} we have $\omega \big|_S =
\pi^{\ast} (\omega_0)_{(\tau)}$ and so for all $\phi \in T_pS$,
\[(d\pi_p(\phi))([f|_{\tau}]) = \phi([\tilde{f}|_S]) =
\omega(\Xi_{\tilde{f}}(p), \phi) =
\pi^*(\omega_0)_{({\tau})}(\Xi_{\tilde{f}}(p), \phi) =
(\omega_0)_{({\tau})}(d\pi_p (\Xi_{\tilde{f}}(p)), d\pi_p(\phi)).\]
It is shown in the proof of Theorem \ref{strat} that the map $\pi: S
\to (M_{\el})_{\tau}$ is a (complex analytic) fibration of type
$G/H$. Thus the differential $d\pi_p : T_pS \to
T_{p_0}((M_{\el})_{\tau})$ is surjective. Therefore
$d\pi_p(\Xi_{\tilde{f}}(p)) = \Xi_{f|_{\tau}}(p_0)$ (the latter
Hamiltonian vector field being defined with respect to
$(\omega_0)_{({\tau})}$).

By definition of the bracket $\{-,-\}_1$ we have $\{f,g\}_1(p_0) =
(\omega_0)_{({\tau})}^i(\Xi_{f|_{\tau}} (p_0), \Xi_{g|_{\tau}}
(p_0))$, and then one calculates
\begin{eqnarray*}\{f,g\}_1(p_0) = (\omega_0)_{({\tau})}(\Xi_{f|_{\tau}} (p_0), \Xi_{g|_{\tau}} (p_0)) =
(\omega_0)_{({\tau})}(d\pi_p(\Xi_{\tilde{f}} (p)),
d\pi_p(\Xi_{\tilde{g}} (p))) \\ =
\pi^*(\omega_0)_{({\tau})}(\Xi_{\tilde{f}}(p), \Xi_{\tilde{g}}(p)) =
\omega(\Xi_{\tilde{f}}(p), \Xi_{\tilde{g}}(p)) = \{\tilde{f},
\tilde{g}\}(p) = \{f,g\}_2(p_0).\end{eqnarray*}
\end{proof}

\subsection{}\label{leavesforlocal}We can now describe the symplectic leaves of
$M_{\el}$. We refer to the Poisson bracket on $\Co(M_{\el})$ from
the proposition above by $\{-,-\}$.
\begin{prop} The orbit type strata $(M_{\el})_{\tau}$ are the symplectic
leaves of $M_{\el}$. There are finitely many symplectic leaves and
they are irreducible locally closed subvarieties.
\end{prop}
\begin{proof}
Let $\tau \in \mathcal{T}$. The set $X = \bigcup_{\nu \leq \tau}
(M_{\el})_{\tau}$ is closed by Proposition \ref{orbittypelocclosed}.
Let $I$ be the ideal of functions vanishing on $X$. Let $x\in X$ so
that $x\in (M_{\el})_{\nu}$ for some $\nu \leq \tau$. For any $f\in
\Co (M_{\el})$ and $i\in I$ we have $\{f,i\}(x) = \{f\big|_{{\nu}}
,i\big|_{{\nu}}\}_{\nu} (x) = \{f\big|_{{\nu}} , 0 \}_{\nu} (x) =
0$, that is, $I$ is a Poisson ideal of $\Co (M_{\el})$. Therefore
the radical of $I$ is a Poisson ideal of $\Co (M_{\el})$ by
\cite[3.3.2]{dix}. Thus $X$ is a Poisson subvariety of $M_{\el}$.

On the other hand, the set $X' = X\setminus (M_{\el})_{\tau}$ is
closed and an identical argument to that of the previous paragraph
shows that it is a Poisson subvariety of $M_{\el}$. Let $p\in
(M_{\el})_{\tau}$ and let $\s$ be the symplectic leaf through $p$ in
$M_{\el}$. If $\s \cap X' \neq \emptyset$ then $\s \subseteq X'$ by
Proposition \ref{algleaves}, but this cannot happen because $p
\notin X'$. Therefore $\s \subseteq X\setminus X' =
(M_{\el})_{\tau}$.

Since $(M_{\el})_{\tau}$ is a symplectic manifold $T_p
(M_{\el})_{\tau}$ is spanned (in the notation of the proposition
above) by the set $\{\{f,-\}_2(p) : f\in \Co (M_{\el})\}$. As we
have shown above, $\{f,-\}_2(p) = \{f,-\}_1(p)$ for all $f\in \Co
(M_{\el})$ and so $T_p (M_{\el})_{\tau}$ is spanned by Hamiltonian
vector fields for the Poisson bracket on $\Co (M_{\el})$. Our choice
of $p\in (M_{\el})_{\tau}$ was arbitrary and $(M_{\el})_{\tau}$ is
connected by Proposition \ref{orbittypelocclosed} so we deduce that
$(M_{\el})_{\tau} \subseteq \s$. Therefore the $(M_{\el})_{\tau}$
are the symplectic leaves of $M_{\el}$ and their properties are a
consequence of Proposition \ref{orbittypelocclosed}.
\end{proof}

\section{Deformed preprojective algebras}\label{RepQ}

\subsection{Representations of quivers}\label{repQdefn}We apply the results of the previous section
to representations of quivers, in particular we prove Theorem
\ref{intro}.

Let $Q$ be a quiver with vertex set $I$, set of arrows $A$, and let
$h, t: A \to I$ denote the head and tail functions respectively. If
$\alpha\in \N^I$, the space of representations of $Q$ of dimension
vector $\alpha$ is
\[
\rp{Q}{\alpha} = \bigoplus_{a\in A} \mat(\alpha_{h(a)}\times
\alpha_{t(a)},\C).
\]
The group
\[
\gp(\alpha) = (\prod_{i\in I} \gl(\alpha_i,\C))/\C^{\times}
\]
acts by conjugation on $\rp{Q}{\alpha}$. It is clear that
$\rp{Q}{\alpha}$ classifies representations of the path algebra $\C
Q$ of dimension vector $\alpha$, and that the $\gp(\alpha)$-orbits
are the isomorphism classes of representations.

\subsection{Moment map}
\label{mom1} Given $\alpha\in\mathbb{N}^I$, $\rp{Q}{\alpha}$ is an
affine space. Let $\overline{Q}$ denotes the double of $Q$ obtained
from $Q$ be adjoining a reverse arrow $a^{\ast} : j\rightarrow i$
for each arrow $a: i\rightarrow j$ in $A$. We denote the set of
arrows in $\overline{Q}$ by $\overline{A}$. Through the trace
pairing,
\[\rp{\overline{Q}}{\alpha} = \bigoplus_{a\in A}
\mat(\alpha_{h(a)}\times \alpha_{t(a)},\C) \oplus \bigoplus_{a\in A}
\mat(\alpha_{h(a^*)}\times \alpha_{t(a^*)},\C)\] is the cotangent
bundle of $\rp{Q}{\alpha}$. 
There is a canonical symplectic form on $\rp{\overline{Q}}{\alpha}$
by
\[\omega(({B}_a, {B}_{a^{\ast}}), ({C}_a, {C}_{a^{\ast}})) = \sum_{a\in A}
-\tr{}{B_{a^{\ast}} C_a} + \tr{}{C_{a^{\ast}} B_a}\] for all
$({B}_a, {B}_{a^*}),({C}_a, {C}_{a^*}) \in
\rp{\overline{Q}}{\alpha}$.

The action of $G(\alpha)$ extends to $\rp{\overline{Q}}{\alpha}$ and
preserves the symplectic form. By Theorem \ref{mw} this action is
Hamiltonian and one can easily verify (see \cite[Page 258]{cb1} for
example) that the corresponding moment map is
\[ \mu_{\alpha}: \rp{\overline{Q}}{\alpha} \longrightarrow
\ed{\alpha}_0, \quad \mu_{\alpha}({B}_{a}, {B}_{a^{\ast}})_i =
\sum_{a\in A, h(a)=i}B_aB_{a^{\ast}} - \sum_{a\in A,
t(a)=i}B_{a^{\ast}}B_a\] where
\[\ed{\alpha}_0 = \{ \theta \in \bigoplus_{i\in I}\mat(\alpha_i \times \alpha_i,\C)
: \sum_{i\in I} \tr (\theta_i) = 0\} \cong (\lie
\gp(\alpha))^{\ast}.\] The isomorphism  $\ed{\alpha}_0 \cong (\lie
\gp(\alpha))^*$ is obtained via the trace pairing and is $\gp
(\alpha)$-equivariant.

Given $\lambda \in \C^I$ with $\lambda \cdot \alpha = 0$, there is a
$\gp(\alpha)$-invariant element of $\ed{\alpha}_0$ whose $i$th
component is $\lambda_i\mathrm{Id}_{\alpha_i}$. The corresponding
Marsden-Weinstein reductions are
\[
\mathcal{N}(\lambda, \alpha) =
\mu_{\alpha}^{-1}(\lambda_i\mathrm{Id}_{\alpha_i})//\gp(\alpha).
\]

\subsection{}\label{preproj}
Given $\lambda\in \C^I$ define the deformed preprojective algebra
\[
\Pi_{\lambda} = \frac{\C\overline{Q}}{(\sum_{a\in A}[a,a^{\ast}] -
\sum_{i\in I}\lambda_ie_i)},
\]
where $[a,a^{\ast}] = aa^{\ast} - a^{\ast}a$. This algebra is
independent of the orientation of $Q$. If $\alpha\in \N^I$ is such
that $\lambda\cdot \alpha=0$ then the representations of
$\Pi_{\lambda}$ of dimension $\alpha$ can be
$\gp(\alpha)$-equivariantly identified with
$\mu_{\alpha}^{-1}(\lambda)$. By a result of Artin, \cite[Section
12]{artin}, $\mathcal{N}(\lambda,\alpha)$ classifies the isomorphism
classes of semisimple representations of $\Pi_{\lambda}$ of
dimension $\alpha$. It is known if $\lambda\cdot \alpha \neq 0$ then
$\Pi_{\lambda}$ has no representations of dimension $\alpha$,
\cite[Theorem 1.2]{cb1}.

\subsection{Stratifying
$\mathcal{N}\left(\lambda,\alpha\right)$ by representation
type}\label{reptype}

Recall the variety $\mathcal{N}(\lambda,\alpha)$ classifies
semisimple $\Pi_{\lambda}$-modules. If $M$ is a semisimple
$\Pi_{\lambda}$-module then we can decompose it into its simple
components $M = M^{\oplus k_1}_1 \oplus \dots \oplus M^{\oplus
k_r}_r$ where the $M_t$ are non-isomorphic simples. If $\beta^{(t)}$
is the dimension vector of $M_t$, then we say $M$ has representation
type \[ (k_1,\beta^{(1)}; \ldots ; k_r,\beta^{(r)}),
\]
which is defined up to permutation of the entries $(k_i,
\beta^{(i)})$. Let $\tau = (k_1,\beta^{(1)}; \ldots ;
k_r,\beta^{(r)})$ and let $\mathcal{R}_{\tau}$ be the subset of
$\mathcal{N}\left(\lambda,\alpha\right)$ with representation type
equal (up to permutation) to $\tau$. We shall refer to the
stratification $\mathcal{N}\left(\lambda,\alpha\right) =
\bigsqcup_{\tau; \mathcal{R}_{\tau} \neq \emptyset}
\mathcal{R}_{\tau}$ as \textit{the stratification by representation
type for} $\mathcal{N}\left(\lambda,\alpha\right)$. Recall the
stratification by orbit type from \ref{orbittype}.

\begin{trm}
The stratification by representation type is equal to the
stratification by orbit type, that is, for each $\tau$,
$\mathcal{R}_{\tau} = \mathcal{N}\left(\lambda,\alpha\right)_{\nu}$
for some $\nu \in \mathcal{T}$.
\end{trm}
\begin{proof}
For any quiver $Q$, the proof of \cite[Theorem 2]{lebrupro} shows
that the stratification by representation type of
$\rep(\overline{Q}, \alpha)//\gp(\alpha)$ is equal to the
stratification by orbit type. Since $\gp(\alpha)$ is reductive one
can identify $\mu^{-1}_{\alpha}(\lambda)//\gp(\alpha)$ with a closed
subvariety of $\rep(\overline{Q}, \alpha)//\gp(\alpha)$. To conclude
we note that the subsets $\mathcal{R}_{\tau},
\mathcal{N}\left(\lambda,\alpha\right)_{\nu}$ are obtained by
intersecting $\mu^{-1}_{\alpha}(\lambda)//\gp(\alpha)$ with the
corresponding representation type and orbit type strata,
respectively, of $\rep(\overline{Q}, \alpha)//\\ \gp(\alpha)$.
\end{proof}

\subsection{Hyper-K\"{a}hler structure}\label{localformQ}
Following Nakajima, \cite[Section 3]{nak4} and \cite[Section
2]{nak3}, we describe the hyper-K\"{a}hler structure of
$\rep(\overline{Q}, \alpha)$, see \ref{hyperkahlerdefn}. We consider
the (real) manifold $\rep (\overline{Q}, \alpha)$. Let $I$ denote
the complex structure given by multiplication by $\sqrt{-1}$. We
define a Hermitian inner product on $\rep (\overline{Q}, \alpha)$ as
follows. For each vertex $i$ we give each $\C^{\alpha_i}$ the
standard Hermitian inner product. We get a Hermitian inner product
on $\mat(\alpha_j \times \alpha_i)$ by $(B,C) = \mathrm{tr}
(BC^{\dag})$ for all $B, C \in \mat(\alpha_j \times \alpha_i)$ (we
use $\dag$ to denote the Hermitian conjugate). Extending this to the
whole of $\rep (\overline{Q}, \alpha)$ yields a Hermitian inner
product, $(-,-)$. Taking the real part of $(-,-)$ yields a real
inner product which is invariant under $I$. We get further
$\R$-linear maps by
\begin{equation*}
J({B_a}, {B_{a^*}}) = ({B_{a^*}^{\dag}}, -{B_a^{\dag}}),\\
K = -JI
\end{equation*}
and it is a straightforward check that $I^2 = J^2 = K^2 = IJK = -
\mathrm{Id}$. Thus $\rep(\overline{Q}, \alpha)$ is a
hyper-K\"{a}hler manifold.

In the notation of \ref{hyperkahlerforms} we have a real symplectic
form $\omega_1$ and a complex symplectic form $\omega$. This complex
symplectic form is the same as the one defined in \ref{mom1} above.
Let $K(\alpha) = \prod_{k\in I} U(\alpha_i)/U(1)$, a maximal
connected (real) subgroup of $\gp (\alpha)$. Then $K (\alpha)$ and
$\gp(\alpha)$ act linearly on $\rep(\overline{Q}, \alpha)$ and
preserve $\omega_1$ and $\omega$ respectively. In this situation the
complex moment map $\mu$ from \ref{hyperkahlerforms} equals
$\mu_{\alpha}$ from \ref{mom1}.

\subsection{}


We have established that $\rep (\overline{Q}, \alpha)$ satisfies all
of the hypotheses of Section \ref{strats} so we can apply all of the
results therein. In particular we can prove Theorem \ref{Qleaves}.

\begin{proof}[Proof of Theorem \ref{intro}]
By Proposition \ref{leavesforlocal} the symplectic leaves of
$\mathcal{N}(\lambda,\alpha)$ are the orbit type strata, and by
Theorem \ref{reptype} the orbit type strata equal the representation
type strata.
\end{proof}

\subsection{Roots}\label{start} Let $Q, A, I, h, t$ be as in \ref{repQdefn}. Our
goal is to use Theorem \ref{Qleaves} to describe the symplectic
leaves of certain deformed quotient singularities. Given that the
leaves of $\mathcal{N}(\lambda, \alpha)$ are described by
representation types it is clear that we will need more information
about the dimension vectors of simple representations. This leads us
to the notion of root vectors of $Q$.

Elements of $\mathbb{Z}^I$ are \textit{vectors} and we write $\ep_i$
for the coordinate vector at vertex $i$. We say that a vector
$\alpha$ has \textit{connected support} if the quiver with vertices
$\{i\in I: \alpha_i \neq 0\}$ and arrows $\{a\in A: \alpha_{h(a)},
\alpha_{t(a)}\neq 0\}$ is connected. We partially order
$\mathbb{Z}^I$ via $\alpha \geq \beta$ if $\alpha_i\geq \beta_i$ for
all $i$, and we write $\alpha > \beta$ to mean that $\alpha\geq
\beta$ and $\alpha \neq \beta$. A vector $\alpha$ is
\textit{positive} if $\alpha > 0$, and \textit{negative} if $\alpha
< 0$. For any subset $X \subseteq \Z^I$ we denote the positive and
negative roots in $X$ by $X^+$ and $X^-$ respectively. We shall
denote the standard inner product of $x,y \in \C^I$ by $x \cdot y$.

The \textit{Ringel form} on $\mathbb{Z}^I$ is defined by
\[ \langle \alpha,\beta \rangle = \sum_{i\in I} {\alpha_i\beta_i} -
\sum_{a\in A}{\alpha_{t(a)}\beta_{h(a)}}.\]
Let $(\alpha,\beta) = \langle \alpha,\beta \rangle + \langle
\beta,\alpha\rangle$ be its symmetrisation. Define \[ p(\alpha) = 1
+ \sum_{a\in A}\alpha_{t(a)}\alpha_{h(a)} - \alpha \cdot \alpha =
1-\frac{1}{2}(\alpha,\alpha).\]

The \textit{fundamental region}, $\mathcal{F}$, is the set $0\neq
\alpha\in \N^I$ with connected support and with $(\alpha ,\ep_i)\leq
0$ for every vertex $i$. If $i$ is a loopfree vertex (so
$p(\ep_i)=0$), there is a reflection $s_i:\mathbb{Z}^I\rightarrow
\mathbb{Z}^I$ defined by
\[s_i(\alpha) = \alpha - (\alpha,\ep_i)\ep_i.\]

The \textit{real roots}, $Re$, are the elements of $\mathbb{Z}^I$
which can be obtained from the coordinate vector at a loopfree
vertex by applying a sequence of reflections at loopfree vertices.
The \textit{imaginary roots}, $Im$, are the elements of $\Z^I$ which
can be obtained from $\mathcal{F}\cup -\mathcal{F}$ by a sequence of
reflections at loopfree vertices. For a quiver without loops it is
easy to see that $\alpha \in Re$ implies that $p(\alpha) =0$. The
set of \textit{roots} is $R = Re \cup Im$. For $\lambda \in \C^I$ we
set $R_{\lambda}= \{ \alpha \in R : \lambda \cdot \alpha = 0\}$.

\subsection{Simple dimension vectors}
\label{simp} For a representation type $\tau =(k_1,\beta^{(1)};
\ldots ; k_r,\beta^{(r)})$ (as defined in \ref{reptype}),
Crawley-Boevey proved the following
\begin{prop}\cite[Theorem 1.3]{cb1}\label{locclosed}
If $\mathcal{R}_{\tau} \neq \emptyset$ then $\mathcal{R}_{\tau}$ is
an irreducible locally closed subset of
$\mathcal{N}\left(\lambda,\alpha\right)$ of dimension $\sum_{t=1}^r
2p(\beta^{(t)})$.
\end{prop}

In order to understand $\mathcal{N}(\lambda,\alpha)$ further we
shall need to describe the dimension vectors of simple
$\Pi_{\lambda}$-modules. Let $R_{\lambda}^+$ be the set of positive
roots $\beta$ with $\beta \cdot \lambda = 0$.

\begin{thm}\cite[Corollary 1.4]{cb1}
For $\lambda\in  \C^I$ and $\alpha\in \N^I$ the following are
equivalent
\begin{enumerate}
\item There is a simple representation of $\Pi_{\lambda}$ of
dimension vector $\alpha$; \item $\alpha\in R_{\lambda}^+$ and
$p(\alpha) > \sum_{t=1}^r p(\beta^{t})$ for any decomposition
$\alpha = \beta^{(1)}+\cdots +\beta^{(r)}$ with $r\geq 2$ and
$\beta^{(t)}\in R_{\lambda}^+$.
\end{enumerate}
In this case $\mathcal{N}(\lambda,\alpha)$ is an irreducible variety
of dimension $2p(\alpha)$ and its general element is a simple
representation of $\Pi_{\lambda}$.
\end{thm}
Henceforth we will let $\Sigma_{\lambda}$ denote the set of
dimension vectors of irreducible representations of $\Pi_{\lambda}$,
that is, vectors satisfying condition $(2)$ of the above theorem.

\subsection{Decomposition theorem}
\label{decomp} Let $\lambda\in \C^I$ and let $\N R_{\lambda}^+$
denote the set of sums (including zero) of the elements of the set
$R_{\lambda}^+$. Let $\alpha\in \N R_{\lambda}^+$. Define
\[
|\alpha|_{\lambda} = \max \left\{ \sum_{t=1}^r p(\beta^{(t)}) :
\alpha = \sum_{t=1}^r\beta^{(t)} \text{ with } \beta^{(t)} \in
\Sigma_{\lambda} \text{ for all }t\right\}.
\]

The following theorem of Crawley-Boevey reduces the study of general
$\mathcal{N}(\lambda,\alpha)$ to the case where $\alpha \in
\Sigma_{\lambda}$. For any variety, $X$, let $\sym^n X$ denote the
$n^{th}$ symmetric product of $X$.
\begin{thm}\cite[Theorem 1.1]{cb2}
Let $\lambda \in \C^I$ and $\alpha \in \N R_{\lambda}^+$. Then
\begin{enumerate}
\item There is a unique decomposition $\alpha = \sigma^{(1)} +
\cdots \sigma^{(r)}$ with $\sigma^{(t)}\in \Sigma_{\lambda}$ for all
$t$, such that $|\alpha|_{\lambda} = \sum_{t=1}^r p(\sigma^{(t)})$.
\item Any other decomposition of $\alpha$ as a sum of elements of
$\Sigma_{\lambda}$ is a refinement of this decomposition. \item
Collecting terms and rewriting this decomposition as $\alpha =
\sum_{t=1}^s m_t \sigma^{(t)}$ where $\sigma^{(1)},\ldots
,\sigma^{(s)}$ are distinct and $m_1,\ldots ,m_s$ are positive
integers, we have
\[
\mathcal{N}(\lambda,\alpha) \cong \prod_{t=1}^s
\sym^{m_t}\mathcal{N}(\lambda,\sigma^{(t)}).
\]
In particular $\mathcal{N}(\lambda,\alpha)$ is irreducible of
dimension $2|\alpha|_{\lambda}$.
\end{enumerate}
\end{thm}


\subsection{The smooth locus of $\mathcal{N}(\lambda,\alpha)$}\label{smooth}
Using the above theorem we can determine whether
$\mathcal{N}(\lambda,\alpha)$ is smooth. The theorem below is known:
for $\alpha \in \Sigma_{\lambda}$, Le Bruyn has proved that the
smooth locus of $\mathcal{N}(\lambda,\alpha)$ is the stratum of
representation type $(1, \alpha)$, \cite[Theorem 3.2]{lebruyn}. For
general $\alpha$ with decomposition $\alpha = \sigma^{(1)} + \cdots
\sigma^{(r)}$ as in Theorem \ref{decomp} (1), collecting terms we
have $\mathcal{N}(\lambda,\alpha) \cong \prod_{t=1}^s
\sym^{m_t}\mathcal{N}(\lambda,\sigma^{(t)})$. Then one observes that
for any affine algebraic variety, $X$, with $\mathrm{Dim}\ X \geq
2$, the smooth locus of $\sym^n X$ is the set $\{ [x_1, \dots ,x_n]
: \mathrm{the}\ x_i\ \mathrm{are\ smooth\ points\ in}\ X\
\mathrm{and}\ x_i \neq x_j\ \mathrm{for\ all}\ i\neq j \}$. It
follows that the smooth locus of $\mathcal{N}(\lambda,\alpha)$ is
the stratum of representation type $(1,\sigma^{(1)};\ldots
;1,\sigma^{(r)})$. We give a proof utilising Theorem \ref{intro}.

\begin{thm}
Let $\lambda \in \C^I$ and $\alpha \in \N R_{\lambda}^+$. The smooth
locus of $\mathcal{N}(\lambda,\alpha)$ coincides with the stratum of
representation type $(1,\sigma^{(1)};\ldots ;1,\sigma^{(r)})$, where
$\alpha = \sigma^{(1)} + \cdots \sigma^{(r)}$ is the unique
decomposition from Theorem \ref{decomp} (1).
\end{thm}
\begin{proof}
Let $\tau = (1,\sigma^{(1)};\ldots ;1,\sigma^{(r)})$, then
$\rm{Dim}\ \mathcal{R}_{\tau} = \rm{Dim}\
\mathcal{N}(\lambda,\alpha)$ by Proposition \ref{locclosed} and
Theorem \ref{decomp} (3). Therefore $\overline{\mathcal{R}}_{\tau} =
\mathcal{N}(\lambda,\alpha)$ because both varieties are irreducible.
By Theorem \ref{intro} and \cite[Proposition 3.7]{brgo},
$\mathcal{R}_{\tau}$ is the smooth locus of
$\mathcal{N}(\lambda,\alpha)$.
\end{proof}

In particular, we have established that the smooth locus of
$\mathcal{N}(\lambda,\alpha)$ is always symplectic.

\subsection{}\label{smooth2} For convenience we spell out how to determine whether
$\mathcal{N}(\lambda,\alpha)$ is smooth in terms of roots of $Q$. We
write $\alpha = \sigma^{(1)} + \cdots \sigma^{(r)} = \sum_{i=1}^t
m_i \sigma^{(i)}$, by collecting like terms.

\begin{cor}
$\mathcal{N}(\lambda,\alpha)$ is smooth if and only if $\alpha =
\sum_{i=1}^t m_i \sigma^{(i)}$ is the only possible decomposition of
$\alpha$ as a sum of elements of $\Sigma_{\lambda}$ and for each
$i$, $p(\sigma^{(i)}) > 0$ implies that $m_i = 1$.
\end{cor}
\begin{proof}
By the theorem, $\mathcal{N}(\lambda,\alpha)$ is smooth if and only
if $(1,\sigma^{(1)};\ldots ;1,\sigma^{(r)})$ is the unique
representation type. The result follows from the lemma below.
\end{proof}

\begin{lem}\cite[Page 260]{cb1}
Let $\beta \in \Sigma_{\lambda}$. If $p(\beta) = 0$ then there is a
unique simple representation of $\Pi_{\lambda}$ with dimension
vector $\beta$ (up to isomorphism); if $p(\beta) > 0$ then there are
infinitely many non-isomorphic simple representations with dimension
vector $\beta$.
\end{lem}

\subsection{Symplectic leaves}\label{leafcalc}

We also show how to work out the symplectic leaves of
$\mathcal{N}(\lambda,\alpha)$. These are equal to the representation
type strata by Theorem \ref{Qleaves} and so this task amounts to
finding each of the possible decompositions of $\alpha$ as a sum of
vectors in $\Sigma_{\lambda}$ and listing the representation types
arising out of each decomposition.

For a decomposition
\begin{equation}\label{reptypedecomp}
\alpha = n_1\gamma_1 + \dots n_s \gamma_s + m_1 \beta_1 + \dots +
m_t \beta_t
\end{equation}
where $\gamma_i, \beta_j \in \Sigma_{\lambda}$ for all $i, j$,
$p(\gamma_i)=0$ and $p(\beta_j)>0$ for all $i,j$, the corresponding
representation types of $\mathcal{N}(\lambda,\alpha)$ are labeled by
$t$-tuples of partitions. More precisely, let $P_j$ be the set of
partitions of $m_j$ for each $j$, and for any $\sigma \in P_j$
denote its length by $l(\sigma)$. All the representation types
coming from $(\ref{reptypedecomp})$ are parametrised naturally by
$P_1 \times \dots \times P_t$, thanks to Lemma \ref{smooth2}. For
each $(\sigma_1, \dots , \sigma_t) \in P_1 \times \dots \times P_t$
we denote the corresponding representation type by $\tau_{(\sigma_1,
\dots , \sigma_t)}$. Then $\mathrm{Dim}\
\mathcal{R}_{\tau_{(\sigma_1, \dots , \sigma_t)}} = 2\sum_{i=1}^t
l(\sigma_i)p(\beta_i)$ by Proposition \ref{simp}.

\subsection{Extended Dynkin diagrams}\label{useful} We give a useful
lemma which applies to the quivers we consider in the proof of
Theorem \ref{intro2}. Let $Q$ be an extended Dynkin diagram (see
\cite{casslo}), oriented to have no cycles (this is no restriction
since the deformed preprojective algebras are orientation
independent).

Let $\delta$ be the minimal positive imaginary root. Then $\delta$
is isotropic, that is, $(\delta, \delta)=0$ and furthermore
$(\delta, \ep_k) = 0$ for all $k\in I$. Any vertex for which
$\delta_k=1$ is called an \textit{extending vertex} - extending
vertices always exist. We can relabel the vertex set $I$ so that $0$
is an extending vertex. All real roots satisfy $(\alpha,\alpha) =
1$.

Let $Q'$ be the quiver obtained from $Q$ by adding a vertex $\infty$
and one arrow from $0$ to this vertex. We'll use apostrophes to
denote data associated with $Q'$.

\begin{lem} Let $n\in \N$. Then
\begin{itemize}
\item[(1)] $
p'(\ep_{\infty} + n\delta) = n.$

\item[(2)] If $\ep_{\infty} + n\delta = \beta^{(1)} + \cdots +
\beta^{(r)}$, where the $\beta^{(i)}$ are positive roots, then
$\sum_{t=1}^r p'(\beta^{(t)}) \leq p'(\ep_{\infty} + n\delta)$, with
equality exactly when all but one of the $\beta^{(t)}$ are equal to
$\delta$.
\end{itemize}
\end{lem}
\begin{proof}\mbox{}
\begin{itemize}
\item[(1)] This is a direct calculation, $p'(\ep_{\infty} + n\delta) = 1 - \frac{1}{2}
((\ep_{\infty},\ep_{\infty}) + 2(n\delta,\ep_{\infty}) +
(n\delta,n\delta)) = 1 - \frac{1}{2}(2 - 2n + 0) = n.$
\item[(2)] This is \cite[Lemma 9.2]{cb1}.
\end{itemize}
\end{proof}

\section{Calogero-Moser space}\label{cmsection}

\subsection{}\label{R^ndecomp} We now begin the proof of Theorem
\ref{intro2}, which will occupy this and the subsequent section. We
fix notation for the remainder of the paper. Let $n$ be an integer
greater than $1$. Let $(L,\omega_L)$ be a 2-dimensional symplectic
vector space, and $\Gamma<\mathrm{Sp}(L)= \mathrm{SL}(2,\C)$ a
finite subgroup. Let $x,y\in L$ be a symplectic basis and let $R$ be
the regular representation of $\Gamma$. Consider the set,
$\mathbf{C}$, of all functions $\Gamma \setminus \{1\} \to \C$ which
are constant on conjugacy classes; we identify this set with
$\C^{d}$ where $d+1$ equals the number of conjugacy classes in
$\Gamma$. For $\underline{c} \in \C^{d}$ and $c_1 \in \C$, we shall
write $\mathbf{c}=(c_1, \underline{c})$.

As a $\Gamma$-module, $R^n$ decomposes into a sum of irreducible
representations: \begin{equation} R^n = \sum_{i=0}^{d} S_i\otimes
\C^{n\delta_i}\end{equation} where the $S_i$ are the irreducible
modules of $\Gamma$ and $\delta_i = \mathrm{Dim}\ S_i$. We shall
assume that $S_0$ is the trivial module. It is clear from the above
that one can identify the group $\aut_{\Gamma}(R^n)$ with
$\hat{G}(n\delta) = \prod_{i=0}^{d} \gl(n\delta_i,\C)$, where
$\delta\in \N^{d+1}$ is the vector with $i$th entry $\delta_i$.

\subsection{Calogero-Moser space}\label{cmspace} Denote by
$e_{\Gamma}\in \edo_{\Gamma}(R)$ the projector onto the trivial
representation. Let $\underline{c}\in \mathbf{C}$ and $c\in
\edo_{\Gamma}(R)$ multiplication by the central element
$\sum_{\Gamma\setminus \{1\}} \underline{c}(\gamma) \gamma$. Let
$\mathcal{O}$ be the $\gl(n,\C)$-conjugacy class formed by all
$n\times n$-matrices of the form $P-\id$, where $P$ is a semisimple
rank one matrix such that $\tr(P) = \tr(\id) = n$. Define
\[ X_{\mathbf{c}} = \{ \nabla \in \hmm{\Gamma}{L}{\edo_{\C}(R^n)} : [\nabla (x),
\nabla (y) ] \in \frac{1}{2}c_1 |\Gamma|\mathcal{O}\otimes
e_{\Gamma} + \id \otimes c \} ,\] where $c_1 \in \C$ and
$\frac{1}{2}c_1 |\Gamma|\mathcal{O}\otimes e_{\Gamma} + \id \otimes
c \subseteq \edo(\C^n)\otimes \edo_{\Gamma}(R) = \edo_{\Gamma}(\C^n
\otimes R) = \edo_{\Gamma}(R^n)$.

There is an action of $\ggp(n\delta) = \aut_{\Gamma}(R^n)$ by
basechange: let $g \in \aut_{\Gamma}(R^n)$ and $\nabla \in
X_{\mathbf{c}}$, then $(g\circ \nabla )(x)= g\nabla (x)g^{-1}$ for
all $x \in L$. The action factors through $\C^{\times}$ (where we
identify $\C^{\times}$ with $\{\lambda \cdot \mathrm{Id}: \lambda
\in \C^{\times}\}$) so that $\gp(n\delta) =
\ggp(n\delta)/\C^{\times}$ acts on $X_{\mathbf{c}}$. The variety
$X_{\mathbf{c}}//\gp(n\delta)$ is called \textit{Calogero-Moser
Space for $\Gamma_n$}. There is a discussion of this space in
\cite[$\S 11$]{eg}.

\subsection{} The vector space $\hmm{\Gamma}{L}{\edo_{\C}(R^n)} =
(L^*\otimes \edo_{\C} (R^n))^{\Gamma}$ is symplectic with form
$\omega_L \otimes tr$, where $tr$ is the symmetric bilinear form
$tr(\phi, \psi) = tr(\phi \psi)$ for $\phi, \psi \in
\edo_{\C}(R^n)$. Clearly the action of $\gp(n\delta)$ preserves this
form. We can identify $\edo_{\Gamma}^0(R^n):= \{A\in \edo_{\C}
(R^n): tr_{R^n} A = 0\}$ with $(\lie \gp(n\delta))^*$ via the trace
pairing. One can check that $[\nabla(x), \nabla(y)]:R^n \to R^n$ is
a $\Gamma$-map and that $tr ([\nabla(x), \nabla(y)]A) = \frac{1}{2}
(\omega_L \otimes tr) (A \cdot \nabla, \nabla)$ for all $\nabla \in
\hmm{\Gamma}{L}{\edo_{\C}(R^n)}$ and $A \in \edo_{\Gamma}^0(R^n)$.
Therefore by Theorem \ref{mw} the action of $\gp(n\delta)$ is
Hamiltonian with moment map $ \nabla \mapsto [\nabla(x), \nabla(y)]$
and so $X_{\mathbf{c}}//\gp(n\delta)$ is a Marsden-Weinstein
reduction and in particular is a Poisson variety.

\subsection{The shifting trick}\label{zeta} We immediately note an
equivalent form for $X_\mathbf{c}$. Let $\nabla \in
\hmm{\Gamma}{L}{\edo_{\C}(R^n)}$, and let $x^*, y^*$ be a dual basis
to $x, y$ so that we can write $\nabla = x^*\otimes \phi +
y^*\otimes \psi$ for some $\phi, \psi \in \edo_{\C}(R^n)$. Then
$[\nabla(x), \nabla(y)] = \phi \psi - \psi \phi$. In this way it is
straightforward to check that $[\nabla(x), \nabla(y)]$ is the map
\[
\xymatrix { [\nabla, \nabla] : R^n \ar^(.5){\zeta\otimes Id}[r] &
L\otimes L\otimes R^n\ar^(.55){Id \otimes \nabla}[r] & L\otimes R^n
\ar^(.6){\nabla}[r] & R^n}
\]
where $\zeta$ is the $\Gamma$-map
\[
\zeta : \C \to L\otimes L,\ \zeta(1) = x\otimes y - y\otimes x.
\]
By the adjunction of $\mathrm{Hom}$ and $\otimes$ we have a
$G(n\delta)$-equivariant isomorphism between \newline
$\hmm{\Gamma}{L}{\edo_{\C}(R^n)}$ and $\hmm{\Gamma}{L\otimes
R^n}{R^n}$. Therefore
\[
X_\mathbf{c}//\gp(n\delta) = \{ \nabla \in \hmm{\Gamma}{L\otimes
R^n}{R^n} : [\nabla, \nabla] \in \frac{1}{2}c_1
|\Gamma|\mathcal{O}\otimes e_{\Gamma} + \id \otimes c
\}//\gp(n\delta).
\]

We shall apply the shifting trick of Lemma \ref{shiftingtrick} to
the Calogero-Moser space $X_\mathbf{c}//\gp(n\delta)$. Let $m=
-\frac{1}{2}nc_1|\Gamma|$. If $m \neq 0$ let $U_m$ be the
$\gl(n,\C)$-conjugacy class formed by all $n\times n$-matrices which
are semisimple rank one and whose trace is equal to $m$. Set $U_0 =
0$. Then
\begin{equation}\label{shiftedcmspace}\begin{split} X_\mathbf{c}//\gp(n\delta) = \{
(\nabla, P) \in \hmm{\Gamma}{L\otimes R^n}{R^n}\times&
U_m \otimes e_{\Gamma}: [\nabla, \nabla] + P =\\
&-\frac{1}{2}c_1 |\Gamma|\id \otimes e_{\Gamma} + \id \otimes c
\}//\gp(n\delta).
\end{split}\end{equation}

\subsection{Linearized Calogero-Moser space}\label{hatX_c}

Ultimately we shall show that $X_\mathbf{c}//\gp(n\delta)$ is
isomorphic to one of the Marsden-Weinstein reductions for quivers
described in Section \ref{RepQ}. The Calogero-Moser space is also a
Marsden-Weinstein reduction but is defined over coadjoint orbit
which a larger than just a single point, the latter being the case
for quiver reductions. Therefore we shall give a linearized version
of Calogero-Moser space. Let
\[
\begin{split} \hat{X_\mathbf{c}} = \{ (\nabla, I,J)\in
\hmm{\Gamma}{L\otimes R^n}{R^n} \oplus (R^n)^{\Gamma} \oplus
(({R^n})^{\ast})^{\Gamma} : [\nabla,\nabla] + (I\otimes J)\\
= -\frac{1}{2}c_1|\Gamma|\id\otimes e_{\Gamma} + \id \otimes
c\}\end{split}
\]
with the action of $\ggp(n\delta)$ on $\hat{X_\mathbf{c}}$ by
basechange. There is a natural symplectic form on
$\hmm{\Gamma}{L\otimes R^n}{R^n} \oplus (R^n)^{\Gamma} \oplus
(({R^n})^{\ast})^{\Gamma}$ (cf. \ref{cmspace}) given by $\omega
((\nabla, I, J),(\nabla' ,I' ,J')) = (\omega_L \otimes tr)(\nabla,
\nabla') - J(I') + J'(I) = tr[\nabla, \nabla'] - J(I') + J'(I)$. The
action of $\gp(n\delta)$ preserves this form so by Theorem \ref{mw}
the action is Hamiltonian. One can easily verify that the map $\rho:
(\nabla, I, J) \mapsto [\nabla,\nabla] + (I\otimes J)$ is the moment
map for this action. Therefore the space
$\hat{X}_\mathbf{c}//\ggp(n\delta)$ is a Marsden-Weinstein
reduction.

\subsection{} The vector space $(R^n)^{\Gamma}$ is equal to $S_0^n = \C^n$.
We form the symplectic vector space $(R^n)^{\Gamma} \oplus
(({R^n})^{\ast})^{\Gamma}$ with the obvious action of $\gl(n,\C)$
preserving this form. The corresponding moment map $\nu :
(R^n)^{\Gamma} \oplus (({R^n})^{\ast})^{\Gamma} \to \mathfrak{gl}_n$
given by $(I,J) \mapsto I\otimes J$ is $\gl(n,\C)$-equivariant and
Poisson by Proposition \ref{moment}. Here, as usual, we identify
$\mathfrak{gl}_n$ with its dual via the trace pairing. Recall that
$m = -\frac{1}{2}nc_1|\Gamma|$.

\begin{prop}\label{cmiso}
There is an isomorphism of Poisson varieties
$\hat{X}_\mathbf{c}//\ggp(n\delta) \cong
X_\mathbf{c}//\gp(n\delta)$.
\end{prop}

\begin{proof}
Suppose first that $c_1 \neq 0$. Taking traces of the defining
equation of $\hat{X}_{\mathbf{c}}$,
\begin{equation}\label{linCMeqn}
[\nabla,\nabla] + I\otimes J = -\frac{1}{2}c_1|\Gamma|\id\otimes
e_{\Gamma} + \id \otimes {c},\end{equation} yields
$\tr_{R^n}(I\otimes J) = \tr_{R^n}(-\frac{1}{2}c_1 |\Gamma|
\id\otimes e_{\Gamma}) + \tr_{R^n}(\id\otimes {c}) = m$, since ${c}$
is traceless on the regular representation. Therefore $I\otimes J
\in \mathfrak{gl}_n$ is a rank one matrix with trace $m$ (which
implies that it is semisimple).

Recall the notation and results from Example \ref{poisson}. In
particular let $W = \C^n \oplus (\C^n)^* = (R^n)^{\Gamma} \oplus
(({R^n})^{\ast})^{\Gamma}$, and let $\mu : W \to \C$ be the moment
map for the $\C^{\times}$-action. Let $W_m$ be the fibre
$\mu^{-1}(m)$. By the previous paragraph
$\rho^{-1}(-\frac{1}{2}c_1|\Gamma|\id\otimes e_{\Gamma} + \id
\otimes \mathbf{c}) \subseteq \hmm{\Gamma}{L\otimes R^n}{R^n} \times
W_m$. The map $\rho$ is constant on $\C^{\times}$-orbits so induces
a map $\overline{\rho}:\hmm{\Gamma}{L\otimes R^n}{R^n} \times
W_m//\C^{\times} \to \ed{n\delta}$ such that
\begin{equation}\label{rhodiagram} \xymatrix{ \hmm{\Gamma}{L\otimes
R^n}{R^n} \times W_m \ar[dr]^{\rho} \ar[d]_{\pi_{\C^{\times}}} & \\
\hmm{\Gamma}{L\otimes R^n}{R^n} \times W_m//\C^{\times} \ar[r]_{\ \
\ \ \ \ \ \ \ \ \ \ \ \overline{\rho}} & \ed{n\delta}}\end{equation}
commutes. Here $\pi_{\C^{\times}}$ is the orbit map for the
$\C^{\times}$ action. Now, $\rho$ is $\ggp (n\delta)$-equivariant
and so is $\pi_{\C^{\times}}$. Therefore the fact that
$\pi_{\C^{\times}}$ is surjective and $\overline{\rho} \circ
\pi_{\C^{\times}} = \rho$ imply that $\overline{\rho}$ is
$\ggp(n\delta)$-equivariant also.

On the other hand, recall the isomorphism $t: W_m//\C^{\times} \to
U_m$ from Example \ref{poisson}; this is a Poisson $\gl
(n,\C)$-equivariant isomorphism. We claim that the following diagram
commutes
\[ \xymatrix{ \hmm{\Gamma}{L\otimes R^n}{R^n} \times W_m//\C^{\times}
\ar[d]_{\mathrm{Id} \times t} \ar[drr]^{\overline{\rho}} & & \\
\hmm{\Gamma}{L\otimes R^n}{R^n} \times U_m \ar[rr]_{\ \ \ \ \ \ \ \
\ \ [-,-] + \iota}& & \ed{n\delta},}\]where $\iota$ is the embedding
$U_m \hookrightarrow U_m \otimes e_{\Gamma} \subset \ed{n\delta}$
and $[-,-] + \iota$ is the map taking $(\nabla, M) \in
\hmm{\Gamma}{L\otimes R^n}{R^n} \times U_m$ to $[\nabla,\nabla] +
\iota(M)$. It follows from the surjectivity of $\pi_{\C^{\times}}$
and (\ref{rhodiagram}) that it suffices to show that $\rho (\nabla,
I, J) = ([-,-] + \iota) \circ (\mathrm{Id} \times t) \circ
\pi_{\C^{\times}} (\nabla, I, J)$ for all $(\nabla,I,J)\in
\hmm{\Gamma}{L\otimes R^n}{R^n} \times W_m$. However, $(\mathrm{Id}
\times t) \circ \pi_{\C^{\times}} (\nabla, I, J) = (\nabla, I
\otimes J)$ which implies that $([-,-]+ \iota) \circ (\mathrm{Id}
\times t) \circ \pi_{\C^{\times}} (\nabla, I, J) = ([-,-]+
\iota)(\nabla, I \otimes J) = [\nabla, \nabla] + I \otimes J = \rho
(\nabla, I, J)$ as required. Therefore $\overline{\rho}$ is a moment
map for the $\ggp(n\delta)$ action by Lemma \ref{shiftingtrick}.
Thus by (\ref{shiftedcmspace}) the map $\mathrm{Id} \times t$
induces a Poisson isomorphism
$\overline{\rho}^{-1}(-\frac{1}{2}c_1|\Gamma|\id\otimes e_{\Gamma} +
\id \otimes {c})//\ggp(n\delta) \cong X_{\mathbf{c}}//\gp(n\delta)$.
Furthermore the canonical isomorphism
${\rho}^{-1}(-\frac{1}{2}c_1|\Gamma|\id\otimes e_{\Gamma} + \id
\otimes {c})//\ggp(n\delta) \cong
\pi_{\C^{\times}}({\rho}^{-1}(-\frac{1}{2}c_1|\Gamma|\id\otimes
e_{\Gamma} + \id \otimes {c}))//\ggp(n\delta) =
\overline{\rho}^{-1}(-\frac{1}{2}c_1|\Gamma|\id\otimes e_{\Gamma} +
\id \otimes {c})//\ggp(n\delta)$ is Poisson by Proposition
\ref{poisson}.

Now suppose that $c_1 = 0$. Then \[\hat{X}_\mathbf{c} = \{ (\nabla,
I,J)\in \hmm{\Gamma}{L\otimes R^n}{R^n} \oplus (R^n)^{\Gamma} \oplus
(({R^n})^{\ast})^{\Gamma} : [\nabla,\nabla] + (I\otimes J) = \id
\otimes c\}.\] Consider the closed subvariety of
$\hat{X}_\mathbf{c}$,
\[Z = \{ (\nabla, I, J) \in \hat{X}_\mathbf{c}: I = J =0\}.\] It is clear
that projecting onto the first component induces a Poisson
isomorphism $Z//\ggp(n\delta) \cong X_\mathbf{c}//\gp(n\delta)$.
Therefore we have the following diagram
\begin{equation}\label{closedembedding} X_\mathbf{c}//\gp(n\delta) \cong
Z//\ggp(n\delta) \hookrightarrow
\hat{X}_\mathbf{c}//\ggp(n\delta).\end{equation} By Proposition
\ref{phisurj} there is a surjective algebra homomorphism $\phi^* :
\Co(X_\mathbf{c}//\gp(n\delta)) \to Z_c$ where the latter algebra is
the centre of a symplectic reflection algebra for a wreath product
and in particular has Krull dimension $2n$, see \ref{centre}.
Therefore $\mathrm{Dim}\ X_\mathbf{c}//\gp(n\delta) \geq 2n$. On the
other hand Theorem \ref{lincalmwrediso} and Lemma
\ref{quivervarietydimn}, which we shall prove below, imply that
$\hat{X}_\mathbf{c}//\ggp(n\delta)$ is irreducible of dimension
$2n$. Therefore $Z//\ggp(n\delta) =
\hat{X}_\mathbf{c}//\ggp(n\delta)$ and the result follows.
\end{proof}

\subsection{}\label{mckaygraph}\label{lambda} We now make the connection between
Calogero-Moser space for $\Gamma_n$ and the quiver reductions from
Section \ref{RepQ}. We shall use the notation introduced in
\ref{R^ndecomp}. We can associate to $\Gamma$ a graph in the
following way. The \textit{McKay graph} of $\Gamma$ is the graph
with vertex set $I = \{0, \dots ,d\}$ and the number of edges
between $i$ and $j$ is the multiplicity of $S_i$ in $S_j \otimes L$.
This graph is extended of type $\tilde{A}, \tilde{D}\ \mathrm{or}\
\tilde{E}$, so it is simply laced and in particular contains no
double edges, see \cite{casslo} for example.

Let $Q$ be the extended Dynkin diagram corresponding to $\Gamma$
under the McKay correspondence, which is given by choosing an
orientation of the McKay graph. Each vertex corresponds to an
irreducible representation of $\Gamma$ and we choose $\delta \in
\N^I$ so that $\delta_i = \mathrm{dim} S_i$ for all $i\in I$. This
is the minimal imaginary root for $Q$ (see \ref{useful}).

Define a linear map
\begin{align*}\lambda:\mathbb{C}^{d+1} \longrightarrow \mathbb{C}^I, \quad \mathbf{c} = (c_1,
\underline{c}) \mapsto \lambda(c_1, \underline{c})_k =&
-\frac{1}{2}c_1 \mathrm{tr}_{S_k} \sum_{\gamma\in \Gamma}\gamma +
\mathrm{tr}_{S_k}\sum_{\gamma\in
\Gamma\setminus\{1\}}\underline{c}({\gamma})(\gamma)\\ =&
-\frac{1}{2}c_1|\Gamma| \mathrm{tr}_{S_k} e_{\Gamma} +
\mathrm{tr}_{S_k} c.\end{align*}

Now consider the quiver $Q'$ defined by adding a vertex $\infty$ to
the quiver $Q$ and adding one arrow from the vertex $\infty$ to the
vertex $0$. Let $I' = I \cup \{\infty\}$ and let
$\lambda'(\mathbf{c}) = (-\lambda(\mathbf{c})\cdot n\delta,
\lambda(\mathbf{c})) \in \C^{I'}$. Then we have one of the
Marsden-Weinstein reductions from Section \ref{RepQ} defined on the
quiver $Q'$:
\begin{equation*}
\begin{split}
\mathcal{N}\left(\lambda'(\mathbf{c}),\ep_{\infty}+n\delta\right)=
\{ (x,i,j) \in \bigoplus_{a\in \overline{A}}
\mat(n\delta_{h(a)}\times n\delta_{t(a)},\C) \oplus \mat(n\delta_0
\times 1, \C)\\ \oplus \mat(1 \times n\delta_0, \C):
\mu_{n\delta}(x) + ij - ji =
\lambda'(\mathbf{c})\}//\gp(\ep_{\infty}+n\delta).
\end{split}
\end{equation*}

\begin{rem}
As noted in \cite[Section 1]{cb1}
$\mathcal{N}\left(\lambda'(\mathbf{c}),\ep_{\infty}+n\delta\right)$
can be described as one of the quiver varieties defined by Nakajima,
see \cite{nak3}.
\end{rem}

\subsection{}\label{quivervarietydimn} We can calculate the
dimension of this Marsden-Weinstein reduction.
\begin{lem}\label{dimcmspace}
$\mathcal{N}\left(\lambda'(\mathbf{c}),\ep_{\infty}+n\delta\right)$
is an irreducible variety of dimension $2n$
\end{lem}
\begin{proof}
By Theorem \ref{decomp} (3) we know that this variety will be
irreducible and have dimension $2n$ as long as
$\ep_{\infty}+n\delta$ is a sum of elements from
$R_{\lambda'(\mathbf{c})}^+$ and that
$|\ep_{\infty}+n\delta|_{\lambda'(\mathbf{c})} = n$. Let $\beta =
\ep_{\infty}+n\delta$. One sees that $\beta$ is in fact a root:
$(\ep_i , \beta) = 0$ for all $i \notin \{0, \infty\}$, $(\ep_0 ,
\beta ) = (\ep_0 , \ep_{\infty}) = -1$ and $(\ep_{\infty} ,\beta) =
2 - n\delta_0 \leq 0$. By the definition of
$|\beta|_{\lambda'(\mathbf{c})}$ and Lemma \ref{useful} (2),
$|\beta|_{\lambda'(\mathbf{c})} = p'(\beta)$ and, as noted in Lemma
\ref{useful} (1), $p'(\beta) = n$.
\end{proof}

%
%

\subsection{} We state some results we shall need. Recall the map $\zeta$ from
\ref{zeta}; dually, the form $\omega_L$ defines a $\Gamma$-map
$\omega_L : L\otimes L \to \C$. By \cite[Lemma 3.1]{cbh} if $M$ and
$N$ are $\C \Gamma$-modules then there is an isomorphism
\[ \hmm{\Gamma}{L \otimes M}{N} \to
 \hmm{\Gamma}{M}{L \otimes N}; \phi \mapsto \phi^{\flat}\] where $\phi^{\flat} =
 (1 \otimes \phi)(\zeta \otimes 1)$. The next result is not
 quite the one proved in \cite{cbh}, but the same proof works in our
 situation.

\begin{prop}\cite[Lemma 3.2]{cbh}\label{tech}
Let $Q$ be a quiver whose underlying graph is a Dynkin graph. Then
one can choose $\{\Theta_a: a\in \overline{A}\}$ such that each
$\Theta_a$ is a basis of $\hmm{\Gamma}{L\otimes S_{t(a)}}{S_{h(a)}}$
and for each $a\in A$
\[ \Theta_a \Theta_{a^*}^{\flat} = \frac{1}{\delta_{h(a)}}1_{S_{h(a)}}\
\mathrm{and}\ \Theta_{a^*} \Theta_a^{\flat} =
-\frac{1}{\delta_{t(a)}} 1_{S_{t(a)}}.\]
\end{prop}

\subsection{} The idea of our proof of the theorem below is based on
\cite[Section 3]{kuz}.

\begin{thm}\label{lincalmwrediso} There is a Poisson isomorphism
$\hat{X}_\mathbf{c}//\ggp(n\delta) \cong
\mathcal{N}\left(\lambda'(\mathbf{c}),\ep_{\infty}+n\delta\right)$.
\end{thm}
\begin{proof} There is an isomorphism $\gp(\ep_{\infty}+n\delta) \cong
\ggp(n\delta)$ given by $(1,g_i)/\C^{\times} \mapsto (g_i)$ and we
shall identify these two groups in what follows. We can identify
$(R^n)^{\Gamma}$ and $((R^n)^*)^{\Gamma}$ with $\mat(n\delta_0\times
1, \C)$ and $\mat(1\times n\delta_0,\C)$ respectively. In this way
and by using the decomposition of $R^n$, (4), we can describe
\[\hmm{\Gamma}{L\otimes R^n}{R^n} \oplus (R^n)^{\Gamma} \oplus
(({R^n})^{\ast})^{\Gamma}\] as
\[\hmm{\Gamma}{
\sum_{i=0}^{d} L\otimes S_i\otimes \C^{n\delta_i}}{\sum_{i=0}^{d}
S_i\otimes \C^{n\delta_i}} \oplus \mat(n\delta_0\times 1, \C) \oplus
\mat(1\times n\delta_0,\C).\] Now the definition of $Q$ from
\ref{mckaygraph} allows us to identify this space with
\begin{equation}\label{decompeqn}\begin{split}\bigoplus_{a\in \overline{A}}
\hmm{\Gamma}{L\otimes S_{t(a)}}{S_{h(a)}}\otimes
\mat(n\delta_{h(a)}\times& n\delta_{t(a)},\C)\\ &\oplus
\mat(n\delta_0\times 1, \C) \oplus \mat(1\times
n\delta_0,\C).\end{split}\end{equation} For each arrow $a\in
\overline{A}$ we choose a $\Theta_a \in \hmm{\Gamma}{L\otimes
S_{t(a)}}{S_{h(a)}}$ as in Proposition \ref{tech}. There is a
$\ggp(n\delta)$-equivariant isomorphism of vector spaces, $\Psi$,
from \[\bigoplus_{a\in \overline{A}} \mat(n\delta_{h(a)}\times
n\delta_{t(a)},\C) \oplus \mat(n\delta_0\times 1, \C) \oplus
\mat(1\times n\delta_0,\C)\] to (\ref{decompeqn}) given by
\begin{eqnarray*}\Psi ((B_a), i, j) = ((\Theta_a \otimes B_a), i,
j).
\end{eqnarray*}

The isomorphism $\gp(\ep_{\infty}+n\delta) \cong \ggp(n\delta)$
gives an isomorphism $\ed{\ep_{\infty} + n\delta}_0 \cong
\ed{n\delta};\\  (M_{\infty}, (M_i)) \mapsto (M_i)$. As noted in
\ref{R^ndecomp} we can identify the groups $\ggp(n\delta)$ and
$\aut_{\Gamma}(R^n)$ and this yields an isomorphism of Lie algebras
$\edo(n\delta) \cong \edo_{\Gamma}{R^n}; (M_i) \mapsto
(\mathrm{Id}_{S_i} \otimes M_i).$ Thus we have an isomorphism
\begin{equation}\label{liealgisos}\ed{\ep_{\infty} + n\delta}_0
\cong \edo_{\Gamma}{R^n}.\end{equation}

We claim that $((B_a),i,j) \in \mu_{\ep_{\infty}+
n\delta}^{-1}(\lambda'(\mathbf{c}))$ if and only if
$\Psi((B_a),i,j)\in \hat{X}_\mathbf{c}$.
Let $((B_a),i,j) \in \mu_{\ep_{\infty}+
n\delta}^{-1}(\lambda'(\mathbf{c}))$. In the isomorphism
(\ref{liealgisos}) the image of $\lambda'(\mathbf{c})$ is
$\bigoplus_{m\in I}{\delta_m}(-\frac{1}{2}c_1 |\Gamma| \mathrm{Id}
\otimes e_{\Gamma} + \mathrm{Id}\otimes {c})|_{S_m \otimes
\C^{n\delta_m}} \in \edo_{\Gamma}{R^n}$. Thus to prove that
$\Psi((B_a),i,j)\in \hat{X}_\mathbf{c}$ we need only verify that for
each $k\in I$,
\[ [(\Theta_a \otimes B_a), (\Theta_a \otimes B_a)]\big|_{ S_k\otimes \C^{n\delta_k}} = \frac{1}{\delta_k} 1_{S_k}
\otimes(\sum_{a\in Q,h(a)=k}  B_a B_{a^{\ast}} - \sum_{a\in Q,
t(a)=k} B_{a^{\ast}} B_a).\] Let $r_k\otimes v_k \in S_k\otimes
\C^{n\delta_k}$, then
\begin{equation}\label{psieqns}\begin{split}&[(\Theta_a \otimes B_a), (\Theta_a \otimes B_a)](r_k\otimes v_k) = (\Theta \otimes B_a)\circ
(\mathrm{Id}\otimes (\Theta_a \otimes B_a))( (x\otimes y
- y\otimes x)\otimes r_k\otimes v_k)\\ &= (\Theta_a \otimes B_a)\big( \sum_{\substack{a\in \overline{A},\\
t(a)=k}} x \otimes \Theta_a (y\otimes r_k)\otimes B_a (v_k) -
\sum_{\substack{a\in \overline{A},\\ t(a)=k}} y \otimes \Theta_a
(x\otimes r_k)\otimes B_a (v_k) \big)\\ &= \sum_{\substack{a\in
\overline{A},\\ t(a)=k}} \sum_{\substack{b\in \overline{A},\\
t(b)=h(a)}} \Theta_b(x\otimes \Theta_a(
y\otimes r_k))\otimes B_b(B_a(v_k))\\ & \hspace{11.5em} - \sum_{\substack{a\in \overline{A},\\
t(a)=k}} \sum_{\substack{b\in \overline{A},\\ t(b)=h(a)}}
\Theta_b(y\otimes \Theta_a(x\otimes r_k))\otimes
B_b(B_a(v_k))\\ &= \sum_{\substack{a\in \overline{A},\\
t(a)=k}} \sum_{\substack{b\in \overline{A},\\ t(b)=h(a)}}
\big[\Theta_b (x\otimes \Theta_a (y \otimes r_k) - y\otimes
\Theta_a(x \otimes r_k))\big]\otimes
B_b(B_a(v_k)).\end{split}\end{equation} Because $[(\Theta_a \otimes
B_a), (\Theta_a \otimes B_a)]$ is a $\Gamma$-map we must have that
the $b$'s appearing in $(\ref{psieqns})$ have $h(b) = k$ and since
$Q$ is a simply laced quiver it follows that $b = a^{\ast}$.
Therefore $(\ref{psieqns})$ equals:
\[\begin{split} &\sum_{\substack{a\in \overline{A},\\ t(a)=k}}
\big[\Theta_{a^{\ast}} (x\otimes \Theta_a (y \otimes r_k) -
y\otimes \Theta_a(x \otimes r_k))\big]\otimes B_{a^{\ast}}(B_a(v_k))\\
&= \sum_{\substack{a\in \overline{A},\\ t(a)=k}}
\big[\Theta_{a^{\ast}} (1 \otimes \Theta_a)(x\otimes y \otimes r_k -
y \otimes x \otimes r_k)\big]\otimes B_{a^{\ast}}(B_a(v_k))\\ &=
\sum_{\substack{a\in \overline{A},\\ t(a)=k}} \Theta_{a^{\ast}}
\Theta_a^{\flat}(r_k)\otimes B_{a^{\ast}}(B_a(v_k)).
\end{split}\] Now by Proposition \ref{tech} this equals
\[
\big( \frac{1}{\delta_k} 1_{S_k}\otimes (\sum_{a\in A, h(a)=k} B_a
B_{a^{\ast}} - \sum_{a\in A, t(a)=k} B_{a^{\ast}} B_a)\big)
(r_k\otimes v_k),
\]
as required. On the other hand, the above calculation shows that if
$\Psi((B_a),i,j)\in \hat{X}_\mathbf{c}$ then for all $k\in I$ we
have $\mu_{\ep_{\infty}+ n\delta}((B_a),i,j)_k =
\lambda(\mathbf{c})_k$. It remains to show that $-ji =
-\lambda(\mathbf{c})\cdot n\delta$. We note that $\mu_{\ep_{\infty}+
n\delta}((B_a),i,j) = \mu_{n\delta}(B_a) + ij - ji$. Now since
$\mu_{n\delta}(B_a)$ is traceless on $\bigoplus_{k\in I}
\C^{n\delta_i}$ we get $ji = \mathrm{tr}\ ij = \mathrm{tr}\
(\mu_{n\delta}(B_a) + ij) = \mathrm{tr}\ \lambda(\mathbf{c}) =
\lambda(\mathbf{c})\cdot n\delta.$

Therefore $\Psi$ induces an isomorphism of $\mu_{\ep_{\infty}+
n\delta}^{-1}(\lambda'(\mathbf{c}))$ to $\hat{X}_\mathbf{c}$, and
since $\Psi$ is equivariant for the group actions we get an induced
isomorphism
\[\tilde{\Psi}:\mathcal{N}\left(\lambda'(\mathbf{c}),\ep_{\infty}+n\delta\right)  \to \hat{X_\mathbf{c}}//\ggp(n\delta).\]

We show that $\tilde{\Psi}$ is Poisson. For all $((B_a), i, j),
((B_a'), i', j')\in \rp{\overline{Q'}}{\ep_{\infty}+n\delta}$,
\begin{equation*} \omega ((B_a), i, j), ((B_a'), i', j')) =
\sum_{k\in I}\big( \sum_{a\in A; t(a)=k} -\tr{}{(B_{a^{\ast}} B'_a)}
+ \sum_{a\in A; h(a)=k} \tr{}{(B'_{a^{\ast}} B_a)}\big) + tr(j'i -
ji')
\end{equation*} which by an easy modification of the calculations
above is equal to
\begin{equation*}\label{cmsympform}\sum_{k\in I} tr_{S_k\otimes \C^{n\delta_k}} [(\Theta \otimes B_a),
(\Theta \otimes B_a')] + tr(i\otimes j' -i'\otimes j).
\end{equation*} By the description of the symplectic form on $\hmm{\Gamma}{L\otimes R^n}{R^n} \oplus (R^n)^{\Gamma} \oplus
(({R^n})^{\ast})^{\Gamma}$ from \ref{hatX_c} it follows that $\Psi$
intertwines symplectic forms. Thus by Proposition \ref{poisson}
$\tilde{\Psi}$ is a Poisson isomorphism.
\end{proof}

\section{Symplectic reflection algebras}\label{srasection}

\subsection{Definition}\label{sra} Let $(L, \omega_L)$ be as in
\ref{R^ndecomp}. The wreath product $\Gamma_n = S_n\ltimes \Gamma^n$
acts on $V:= L^{\oplus n} = \C^n\otimes L$, preserving the
symplectic form $\omega = \omega_L^n$.

We say that $\gamma \in \Gamma_n$ is a \textit{symplectic
reflection} if $\mathrm{Dim}(1 - \gamma)(V) = 2$. Let $\mathit{S}$
denote the set of all symplectic reflections. Let
$\mathbf{c}:\mathit{S} \to \C; \gamma \mapsto c_{\gamma}$ be
constant on conjugacy classes. Given $\gamma \in \mathit{S}$ define
the form $\omega_{\gamma}$ on $V$ to have radical $\mathrm{ker}(1 -
\gamma)$ and to be the restriction of $\omega$ on $(1 - \gamma)(V)$.

Let $TV$ be the tensor algebra of $V$. The \textit{symplectic
reflection algebra} $H_{\mathbf{c}}$ is the $\C$-algebra defined as
the quotient of the skew group ring $TV*\Gamma_n$ by the relations
\[ v\otimes w - w\otimes v = \sum_{\gamma \in \mathit{S}}
c_{\gamma}\omega_{\gamma}(v,w)\gamma, \] for all $v,w \in V$.

Usually, symplectic reflection algebras depend on a further
parameter $t\in \C$. The definition above is the $t=0$ case.

\subsection{The geometry of the centres}\label{centre} Let $Z_{\mathbf{c}}$ denote
the centre of $H_{\mathbf{c}}$. By the PBW theorem \cite[Theorem
1.3]{eg}, $\mathrm{GKdim}H_{\mathbf{c}} = 2n$. $H_\mathbf{c}$ is a
finitely generated module over $Z_\mathbf{c}$, \cite[Theorem
3.1]{eg}, and it follows that $\mathrm{GKdim}Z_{\mathbf{c}} = 2n$.
By the Artin-Tate Lemma, \cite[Lemma 13.9.10]{mcrob}, $Z_\mathbf{c}$
is a finitely generated algebra over $\C$; also, it follows from
\cite[Theorem 3.3]{eg} that $Z_{\mathbf{c}}$ is a domain. Thus, $\mx
Z_{\mathbf{c}}$ is an irreducible variety of dimension $2n$, and
since $\mx Z_{\mathbf{0}} = V/\Gamma_n$ the varieties $\mx
Z_{\mathbf{c}}$ form a flat family of deformations of the symplectic
quotient singularity $V/\Gamma_n$, \cite[Theorem 3.1]{eg}. The case
$n=1$ is the Kleinian singularity case and that is dealt with in
detail in \cite{casslo}.

By \cite[Section 15]{eg} the algebras $Z_{\mathbf{c}}$ have a
Poisson bracket which is a deformation of the bracket on $\Co
(V/\Gamma_n)$. Our interest is in determining for which values of
$\mathbf{c}$ the varieties $\mx Z_{\mathbf{c}}$ are smooth, and in
the singular cases calculating the symplectic leaves, which have a
leading role in the representation theory of $H_{\mathbf{c}}$.


\subsection{Symplectic reflections}\label{wreath}

We list the symplectic reflections of $\Gamma_n$ for $n>1$. Given
$\gamma\in \Gamma$, we write $\gamma_i \in \Gamma_n$ for the element
$\gamma$ regarded as an element of the $i$-th factor $\Gamma$. Let
$s_{ij}\in S_n$ denote the transposition swapping $i$ and $j$. Then
the group $\Gamma_n$ is generated by the symplectic reflections
$s_{ij}$ and $\gamma_i$. The conjugacy classes of symplectic
reflections in $\Gamma_n<\mathrm{Sp}(V)$ are known to be of two
types:
\begin{enumerate}
\item The set
$\{s_{ij}\cdot\gamma_i\cdot\gamma_j^{-1}: i,j\in{[}1,n{]}, i\neq j,
\gamma\in\Gamma\}$ forms a single $\Gamma_n$ conjugacy class; \item
The set $\{\gamma_i: i\in {[}1,n{]}, \gamma\in \mathcal{C}\}$ forms
one $\Gamma_n$ conjugacy class for any given conjugacy class,
$\mathcal{C}$, in $\Gamma\setminus\{1\}$.
\end{enumerate}

In particular, we can identify the set of class functions $S \to \C$
with $\C^{d+1}$ where $d+1$ is the number of conjugacy classes in
$\Gamma$. Therefore the flat family of (centres of) symplectic
reflection algebras is parametrised by elements $\mathbf{c} \in
\C^{d+1}$. We will assume that $\mathbf{c}= (c_1, \underline{c})$
where $c_1 \in \C$ corresponds to the conjugacy class $(1)$ above
and $\underline{c}\in \C^{d}$ corresponds to class functions $\Gamma
\setminus \{ 1\} \to \C$ cf. \ref{lambda}.

\subsection{} Let $\mx Z_{\mathbf{c}}$ be the variety associated to
$H_{\mathbf{c}}$ as in \ref{centre}. The theorem we now prove is
based very closely on \cite[Theorem 11.16]{eg}. We make use of Lemma
11.15 in \cite{eg} as in the earlier proof, but a check shows that
$c_1$ should be replaced by $\frac{1}{2}c_1$. This accounts for the
discrepancy between our definition of Calogero-Moser space,
$X_\mathbf{c}//\gp(n\delta)$, given in \ref{cmspace} and that given
in \cite[Definition 11.5]{eg}.

\begin{trm}\label{sra-cmiso} $\mathrm{Max}Z_\mathbf{c}$ and
$X_{\mathbf{c}}//\gp(n\delta)$ are isomorphic as Poisson varieties
(up to nonzero scalar multiple).
\end{trm}

\begin{proof}

We first show that $X_{c}//\gp(n\delta)$ is isomorphic to
$\mathrm{Max}Z_\mathbf{c}$. Our argument follows that of
\cite[Theorem 11.16]{eg} which proves that these varieties are
isomorphic for generic values of $\mathbf{c}$. We show that the
generic hypothesis can be removed.

Let $\rep_{\Gamma_n}(H_\mathbf{c})$ be the variety of
$H_\mathbf{c}$-modules which are isomorphic to the regular
representation of $\Gamma_n$. Let $e\in \C \Gamma_n$ be the
symmetrising idempotent. Given $M\in \rep_{\Gamma_n}(H_\mathbf{c})$,
$eM$, the space of $\Gamma_n$ fixed points, is a one dimensional
subspace of $M$ fixed under the action of $Z_\mathbf{c}$. This
induces a morphism of algebraic varieties
\[
\pi : \rep_{\Gamma_n}(H_\mathbf{c}) \longrightarrow \mx
Z_\mathbf{c};\ M \mapsto \mathrm{ann}_{Z_\mathbf{c}} (eM)
\]

The group $\aut_{\Gamma_n}(\C\Gamma_n)$ acts naturally on
$\rep_{\Gamma_n}(H_\mathbf{c})$. By \cite[Theorem 3.7]{eg} there is
a unique irreducible component $\rep_{\Gamma_n}^o(H_\mathbf{c})$ of
the variety $\rep_{\Gamma_n}(H_\mathbf{c})$ whose image is dense in
$\mx Z_\mathbf{c}$ and whose generic point is a simple
$H_\mathbf{c}$-module. Moreover, there is an isomorphism
\[
\pi^{\ast}: Z_\mathbf{c} \longrightarrow \Co (
\rep_{\Gamma_n}^o(H_\mathbf{c}))^{\aut_{\Gamma_n}(C\Gamma_n)}.
\]

On the other hand, thanks to \cite[pp.311]{eg}, there is a morphism
\[
\psi: \rep_{\Gamma_n}(H_\mathbf{c}) \longrightarrow X_{\mathbf{c}}.
\]
This morphism is obtained as follows. Let $S_{n-1}\ltimes
\Gamma^{n-1}= \Gamma_{n-1}<\Gamma_n$, where $S_{n-1}$ fixes the
label $1\in \{1,\ldots ,n\}$. Then, for $v\in L$, the element
$v_1=(v,0,\ldots ,0)\in L^{\oplus n}\subset H_\mathbf{c}$ commutes
with $\Gamma_{n-1}$. Let $M\in \rep_{\Gamma_n}(H_\mathbf{c})$, and
consider the $\Gamma_{n-1}$ fixed points $M^{\Gamma_{n-1}}\subseteq
M$. We have a $\Gamma$-equivariant map
\[ \nabla_M : L \to \edo_{\C}(M^{\Gamma_{n-1}});\ v \mapsto
v_1|_{M^{\Gamma_{n-1}}}.
\]
Now $\psi$ sends $M$ to $\nabla_M$. By \cite[Lemma 11.15]{eg} this
is a well-defined map from $\rep_{\Gamma_n}(H_\mathbf{c})$ to
$X_\mathbf{c}$.

The morphism $\psi$ intertwines the $\aut_{\Gamma_n}(\C\Gamma_n)$
and $\gp(n\delta)$ actions, and hence induces a map
\begin{equation}\label{phidefn} \phi^{\ast} : \Co
(X_{\mathbf{c}}//\gp(n\delta)) \longrightarrow \Co
(\rep_{\Gamma_n}^o(H_\mathbf{c}))^{\aut_{\Gamma_n}(\C\Gamma_n)}
\cong Z_\mathbf{c} .
\end{equation}
By Proposition \ref{phisurj} below $\phi^{\ast}$ is surjective.
However, the proof given in \cite{eg} that $\phi^{\ast}$ is
injective requires the generic hypothesis on $\mathbf{c}$. We
circumvent this in the following way: by combining Lemma
\ref{quivervarietydimn}, Proposition \ref{cmiso} and Theorem
\ref{lincalmwrediso} we have that $X_{\mathbf{c}}//\gp(n\delta)$ is
an irreducible variety of dimension $2n$. The variety
$\mathrm{Max}Z_\mathbf{c}$ is also irreducible of dimension $2n$ and
therefore the surjective morphism $\phi^{\ast}$ must be injective
also. Therefore $X_{\mathbf{c}}//\gp(n\delta) \cong
\mathrm{Max}Z_\mathbf{c}$.

To see that this isomorphism is Poisson (up to nonzero scalar
multiple) we simply note that the penultimate two paragraphs of
\cite[Lemma 11.18]{eg} apply word-for-word: the map $\phi^*$ is a
filtration preserving isomorphism, and the Poisson bracket on
$X_{\mathbf{c}}//\gp(n\delta)$ has filtration degree $\leq 2$ (in
the sense of \cite[Page 267]{eg}). Therefore, by \cite[Lemma
2.26]{eg}, $\phi^*$ is Poisson up to nonzero scalar multiple.
\end{proof}

\subsection{}\label{phisurj} There is a result used in the proof of Theorem
\ref{sra-cmiso} which we require in the proof of Proposition
\ref{cmiso}. The proof given in \cite{eg} can be applied verbatim to
our situation.

\begin{prop}\cite[Theorem 11.16]{eg}
The algebra homomorphism $\phi^*$ from (\ref{phidefn}) is a
surjective map.
\end{prop}

One consequence of Theorem \ref{sra-cmiso} is that the isomorphism
$\mathrm{Max}Z_\mathbf{c} \cong X_{c}//\gp(n\delta)$ maps symplectic
leaves to symplectic leaves, \cite[Lemma 3.4]{mma}.

\section{Examples}\label{examples}
\subsection{The symmetric group}\label{subgroups}
We calculate some examples of Theorem \ref{intro2}, working out for
which values of $\mathbf{c}$ the variety $\mx Z_{\mathbf{c}}$ is
smooth. By \cite[Proposition 3.7 and Theorem 7.8]{brgo}, $\mx
Z_{\mathbf{c}}$ is smooth if and only if it is symplectic and so
there will only be nontrivial leaves when $\mx Z_{\mathbf{c}}$ is
singular. In these cases we work out the number and dimension of the
symplectic leaves. We refer the reader to the notation introduced in
\ref{sra}.

The simplest case to consider is when the group $\Gamma$ is trivial,
that is, where the there is only the symmetric group $S_n$ acting by
its permutation representation on $(\C^2)^n$. Then the parameter
$\mathbf{c}$ is simply a complex number. For generic values of
$\mathbf{c}$ the variety $\mx Z_{\mathbf{c}}$ is smooth by
\cite[Corollary 1.14]{eg}. It is clear from the relations that
$H_{t\mathbf{c}} \cong H_{\mathbf{c}}$ for all $t\neq 0$, thus $\mx
Z_{t\mathbf{c}} \cong \mx Z_{\mathbf{c}}$ and we deduce that $\mx
Z_{\mathbf{c}}$ is smooth for all nonzero $\mathbf{c}$. It is
possible to show that in this case $\mx Z_{\mathbf{c}}$ is
isomorphic to $\mathbb{A}^2 \times \mx \mathfrak Z_{\mathbf{c}}$,
where $\mathfrak Z_{\mathbf{c}}$ is the centre of the rational
Cherednik algebra (at $t=0$) of type $\mathbf{A}_{n-1}$, see
\cite[Section 4]{eg}.

\subsection{The hyperoctrahedral group}
The next simplest case is $\Gamma_n = S_n \ltimes (\Z_2 \times \dots
\times \Z_2)$ where we identify $\Z_2$ with $<\gamma = -
\mathrm{Id}_{\C^2}>$. We consider the family of varieties $\mx
Z_{\mathbf{c}}$ where the deformation parameter is given by two
variables $\mathbf{c} = (c_1,c_{\gamma})$, corresponding to the
conjugacy classes of symplectic reflections in $\Gamma_n$
(\ref{wreath}). Here $\mx Z_{\mathbf{c}}$ is isomorphic to $\mx
\mathfrak Z_{\mathbf{c}}$ where $\mathfrak Z_{\mathbf{c}}$ is the
centre of the rational Cherednik algebra of type $\mathbf{B}_n$.

Theorem \ref{intro2} gives us $\mx Z_{ \mathbf{c}} \cong
\mathcal{N}\left(\lambda'(\mathbf{c}),\ep_{\infty}+n\delta\right)$
and this isomorphism identifies symplectic leaves. The latter is the
Marsden-Weinstein reduction associated to the quiver
\[\mathrm{Q':}\ \xymatrix{ \infty \ar[r] & 0 \ar@/^/[r] & 1 \ar@/^/[l]}
\]
with dimension vector $\ep_{\infty}+n\delta = (1,n,n)$ and at
parameter $\lambda'(\mathbf{c}) = (nc_1, -c_1 + c_{\gamma},
-c_{\gamma})$ (see \ref{mckaygraph}). We note that it is well known
that the variety $\mx Z_{\mathbf{c}}$ is singular when $\mathbf{c}$
equals zero and its symplectic leaves have been calculated,
\cite[Proposition 7.4]{brgo}. We omit the proof of the following
theorem, details of which can be found in the author's PhD thesis.

\begin{thm} Let $\mx Z_{\mathbf{c}}$ be as above and assume that
$\mathbf{c} \neq (0,0)$. Then
\begin{enumerate} \item $\mx Z_{\mathbf{c}}$ is singular if and only
if $c_1 = 0$ or $c_{\gamma} = \pm mc_1$ for some integer, $m$, such
that $0\leq m \leq n-1$.
\item If $c_1 = 0$ then the symplectic leaves of $\mx
Z_{\mathbf{c}}$ are parametrised by $P_n$, the set of partitions of
$n$. For each $\sigma \in P_n$, the corresponding leaf,
$\mathcal{S}_{\sigma}$ has dimension $2l(\sigma)$, where $l(\sigma)$
is the length of $\sigma$.
\item If $c_{\gamma} = \pm mc_1$ for some $0\leq m \leq n-1$ then
the leaves are parametrised by the set $S = \{k \in \Z_{\geq 0} : km
+ k^2 \leq n \}$. For $k \in S$ the corresponding leaf,
$\mathcal{S}_k$ has dimension $2(n - km -k^2)$.
\end{enumerate}
\end{thm}

\begin{rem}
The behaviour of $\mx Z_{\mathbf{c}}$ when $c_1 = 0$ is always
analogous to that occurring in the above theorem. Denote by
$H_{\underline{c}}$ the symplectic reflection algebra associated to
the pair $(L,\omega_L)$ acted on by the group $\Gamma$ (see
\ref{R^ndecomp}), with centre $Z_{\underline{c}}$. When $c_1 = 0$ it
is easy to calculate that $H_{\mathbf{c}}$, the algebra defined on
the wreath product, $\Gamma_n$, is isomorphic to $\big(\bigotimes_n
H_{\underline{c}} \big) * S_n$. Therefore one can see that $\mx
Z_{\mathbf{c}} \cong \sym^n \mx Z_{\underline{c}}$. This is
identical to the behaviour of symplectic reflection algebras for
wreath products in the case that $t \neq 0$, see \cite{ganginzburg}.
\end{rem}

\textit{Acknowledgements} The author wishes to thank Iain Gordon for
suggesting this topic and contributing many helpful ideas and
comments. Many thanks also go to Ken Brown for his suggestions and
advice. This paper will form part of the author's PhD thesis and he
gratefully acknowledges the support of the University of Glasgow.
Part of the work for this paper was during a workshop at the
University of Edinburgh in June 2005, which was supported by
Leverhulme Research Interchange Grant F/00158/X.

\end{document}